\documentclass{article}
\usepackage{amsmath,amstext,amsthm,amsfonts}
\usepackage[dvips]{graphicx}
\usepackage{amssymb}
\usepackage{xcolor}
\usepackage{environ,hyperref}
\numberwithin{equation}{section}
\newtheorem{teorema}{Theorem}[section]
\newtheorem{proposition}[teorema]{Proposition}
\newtheorem{remark}[teorema]{Remark}

\newtheorem{lema}[teorema]{Lemma}
\newtheorem{definicao}[teorema]{Definition}
\newtheorem{corolario}[teorema]{Corollary}

\newcommand{\EE}{\mathbb{E}}
\newcommand{\VV}{\mathbb{V}}
\newcommand{\RR}{\mathbb{R}}

\newcommand{\ZZ}{\mathbb{Z}}
\newcommand{\PP}{\mathbb{P}}
\newcommand{\NN}{\mathbb{N}}

\newcommand{\cC}{\mathcal{C}}
\newcommand{\cE}{\mathcal{E}}
\newcommand{\Ga}{\Gamma}
\newcommand{\F}{\cal F}
\newcommand{\A}{\mathcal{A}}
\newcommand{\B}{\mathcal{B}}

\newcommand{\p}{\mathfrak{p}}
\newcommand{\q}{\mathfrak{q}}

\theoremstyle{plain}
\begin{document}

\title{{{Dependent Percolation on $\ZZ^2$}}}
\author{ B.N.B. de Lima\footnote{Departamento de Matemática. Universidade Federal de Minas Gerais, Av. Ant\^onio Carlos 6627, CEP 30123-970, Belo Horizonte, MG, Brasil. Email:
\textsf{bnblima@mat.ufmg.br}} , \ V. Sidoravicius\footnote{Shanghai New York University, Pudong New District, Shanghai, China.},  M. E. Vares\footnote{Instituto de Matem\'atica. Universidade Federal do Rio de Janeiro, Av. Athos da Silveira Ramos 149, CEP 21941-909, Rio de Janeiro, RJ, Brasil. Email:
\textsf{eulalia@im.ufrj.br}} } \maketitle
\begin{abstract} {We consider a dependent percolation model on the square lattice $\ZZ^2$. The range of dependence is infinite in vertical and horizontal directions. In this context, we prove the existence of a phase transition. The proof exploits a multi-scale renormalization argument that is defined once the environment configuration is suitably good and, which, together with the main estimate for the induction step, comes from Kesten, Sidoravicius and Vares (To appear in {\em Electronic Journal of Probability}, (2022)). This work was inspired by de Lima (Ph.D.Thesis, \emph{Informes de Matemática. IMPA}, Série C-26/2004) where the simpler case of a deterministic environment was considered. It has various applications, including an alternative proof for the phase transition on the two dimensional random stretched lattice proved by Hoffman ({\em  Comm. Math. Phys.} {\bf 254}, 1-22 (2005)).}
\end{abstract}


\section{Introduction}

Questions involving percolation and the behaviour of growth processes in random environment are very natural and have been object of intense research activity. In this paper we continue an investigation project that was initiated in \cite{Lima} and which consists in the treatment of bond percolation on the square lattice for a class of environments that present infinite range dependencies. This class of problems poses different levels of challenge in their rigorous analysis. One can find motivation back in the study of disordered Ising models, as initiated by McCoy and Wu \cite{MW1,MW2,MW3,MW4},
for which, fixed an inverse temperature in the phase transition region, one may consider lowering the interactions $J_{x,y}$  to a lower value, constant along a random set of horizontal lines, and investigating when the phase transition persists. Another class of examples found motivation in problems that originated in theoretical computer science \cite{Wi, KLSV,BS} or we may simply have the issue of modelling percolation in an environment that exhibits a layered structure, with big difference in porosity depending on such layers. This is for instance the situation with the so called Randomly Stretched Lattice treated by Jonasson, Mossel and Peres \cite{JMP} for $\mathbb{Z}^d$ with $d \ge 3$ and by Hoffman \cite{Ho} for $\mathbb{Z}^2$; and \cite{DHKS}. It is also natural to consider similar problems in the oriented percolation, like \cite{KSV}, or contact process setup, for which we mention \cite{BDS}, \cite{Kl} and \cite{Li}. As in most of the previous articles, our proofs rely on a multi-scale renormalization. Our renormalization method is inspired by \cite{KSV} which we use different ways: for a suitable grouping procedure of the environment configuration and then for the induction step of the estimates.

We now give the precise definition of our percolation model and state our main theorem. The graph $\mathbb{L}^2=(\ZZ_+^2, \EE)$ corresponds to the first quadrant of the square lattice with the set of nearest neighbour bonds   $\EE=\{\langle(x_1,x_2);(y_1,y_2)\rangle\in\ZZ_+^2\times\ZZ_+^2:|x_1-y_1|+|x_2-y_2|=1\}$. We set, for each $i\in\ZZ_+$:

$$\EE_i^H=\{\langle(i,y);(i+1,y)\rangle\in\EE:
y \in \ZZ_+\}$$ and
$$\EE_i^V=\{\langle(x,i);(x,i+1)\rangle\in\EE:
x \in \ZZ_+\},$$ 
i.e. $\EE_i^H$  corresponds to a ladder of horizontal bonds and $\EE_i^V$ corresponds to a ladder of vertical bonds in $\EE$.

The model $\PP_{\delta,p_g,p_b}$ has three parameters. The parameter $\delta$ controls the random environment and sets the infinite range dependency. On a suitable probability space $(\Omega, \mathcal{A}, \PP)$, where $\PP=\PP_{\delta,p_g,p_b}$, we consider a pair of independent Bernoulli sequences  $\xi=(\xi^H, \xi^V)$, with $\xi^H=(\xi^H_i\colon i \in \ZZ_+)$ and $\xi^V=(\xi^V_i\colon i \in \ZZ_+)$, and $\PP(\xi^H_i=1)=\PP(\xi^V_i=1)=\delta=1-\PP(\xi^H_i=0)=\PP(\xi^V_i=0)$. The random variable $\xi^H_i$ determines the corresponding ladder $\EE_i^H$ to be {\em good} or {\em bad} according to $\xi^H_i=0$ or $\xi^H_i=1$, respectively. Similarly for each ladder $\EE_i^V$, according to the value of $\xi^V_i$. 

Given the sequences $\xi=(\xi^H, \xi^V)$, we say that all bonds on the vertical line $\{(i,y): y\in\RR\}$ are \emph{useless} if $\xi^H_{i-1}=\xi^H_{i}=1$, analogously all bonds on the horizontal line $\{(x,i): x\in\RR\}$ are useless if $\xi^V_{i-1}=\xi^V_{i}=1$. In words, a bond $e$ is useless if its four orthogonal neighbour bonds belong to bad ladders. Otherwise, we say that the bond $e$ is useful.

Given the random environment $\xi=(\xi^H, \xi^V)$ and the other two parameters $0<p_b< p_g<1$, we define a percolation model by associating for each bond the state {\em open} or {\em closed}, independently, and with the following conditional probabilities, for $i\in\ZZ$ and $j=H,V$: 
\begin{equation}
\label{def-modelo}
\PP(f\mbox{ is open}|\xi)=\left\{
\begin{array}
[c]{l}
0,\mbox{ if}\ f\mbox{ is useless},\\
p_g,\mbox{ if}\ f\mbox{ is useful},\ f\in \EE^j_i\mbox{ and }\xi^j_i=0,\\
p_b,\mbox{ if}\ f\mbox{ is useful},\ f\in \EE^j_i\mbox{ and }\xi^j_i=1.
\end{array}\right.
\end{equation}

In particular, the status open or closed  of any two useful bonds on a same ladder are dependent no matter how distant they are apart.

We can now state the main result of the paper, giving sufficient conditions on  $p_b,p_g$ and $\delta$ for the origin to percolate with positive probability. Along this paper we use the standard notation in
percolation $\{x\leftrightarrow y\}$  to denote the set of configurations for which the vertex $x$ is connected to the vertex $y$ by a path of open bonds, and $\{x\leftrightarrow \infty\}$ denotes the set of configurations for which $x$ is connected to infinitely many vertices by paths of open bonds; and $p_c:=\tfrac{1}{2}$ is the percolation threshold for the square lattice.

\begin{teorema}\label{main}
%There exists $\widehat{p}<1$ so that for all $p_g>\widehat{p}$ 
For any $p_g>p_c$ and $p_b>0$, we can find $\widehat{\delta}(p_b,p_g)>0$ such that $\PP\{0\leftrightarrow\infty\}>0$ whenever $\delta<\widehat{\delta}(p_b,p_b)$. More precisely, for such a choice of the parameters, 
\begin{equation}
\label{eq-main}
\PP\{0\leftrightarrow\infty| \xi\}>0 \text{ for almost all  } \xi.
\end{equation}
\end{teorema}

The proof of Theorem \ref{main} relies on a multiscale renormalization argument based on a grouping procedure applicable to the environment $\xi$, as introduced in \cite{KSV}.
%The next step is to define the renormalized lattice and finally we will show how to drill this grouping of bad ladders, the details of this drilling can also be seen in \cite{KSV}.
For the renormalized percolation model, we found convenient to use a site-bond version with bad layers being treated as bonds.
%, and start by recalling a few basic facts about this model.

One consequence of the theorem above is an alternative proof of phase transition on the Randomly Stretched Lattice in $d=2$, as showed by Hoffman~\cite{Ho}. 

Another consequence of the main theorem is the existence of a phase transition in a two dimensional Potts model in random environment. This is very similar to the dependent Ising model studied by McCoy and Wu mentioned above. That is, similarly to Theorem~\ref{main}, it guarantees a phase transition for a two dimensional Potts model in a random environment, where there are two values of (ferromagnetic) coupling constants which are chosen at random following the same
pattern used in our percolation model. More precisely: given the parameter $\delta\in(0,1)$ and the positive coupling constants $J_1$ and $J_2$, $J_1<J_2$ we associate for each ladder the state
good (or 0) or bad (or 1), independently, with probability $1-\delta$ and $\delta$, respectively. As before, let us denote by $(\xi_e)_{e\in\EE(\ZZ^2)}\in\{0,1\}^{\EE(\ZZ^2)}$ the configuration or good and bad bonds. Given the
configuration of all good and bad bonds $(\xi_e)_{e\in\EE(\ZZ^2}$, we associate for each bond the coupling constants $J_1$ or $J_2$, if this bond belongs to a bad ladder or a good one, respectively. Let $\mu^s_{\delta ,J_1,J_2}$ be the probability measure for this infinite volume Potts model with pure $s$ boundary condition. Thus, for any $J_1>0$ and $J_2$ large enough, there is a small $\delta=\delta(J_1,J_2)>0$, such that the model exhibits a phase transition; a statement analogous to Theorem \ref{main} follows if we use the lemma below (built on the Fortuin-Kasteleyn measure and ferromagnetic inequalities, and which relates the independent percolation probability measure, $P$, and the Gibbs measure, $\mu_{J}^{s}$, for the Potts model. 

\begin{lema}\label{lemfond}
Define $p(e)$ by
\begin{equation}\label{equpJ}
p(e)=\frac{1-\exp[-2 J(e)]}{1+(q-1)\exp[-2 J(e)]}\,.
\end{equation}
Then the origin's magnetization of the $q$-states
Potts model with couplings $J(e)$ and the probability of percolation of the origin
in the independent percolation process with bond probabilities $p(e)$
are related by the following inequality:
\begin{equation}\label{inequfond}
\mu_{J}^{s}(\sigma_0=s)\geq\frac{1}{q}+\frac{q-1}{q}
P(0\leftrightarrow\infty)\,.
\end{equation}
Where $J(e)$ and $p(e)$ are the coupling constant and the probability of the bond $e$ to be open, for the Potts and Percolation models, respectively.
\end{lema}

The proof of Lemma \ref{lemfond} can be found in \cite{ACCN} or \cite{GHM}.

In Section~\ref{site-bond}, we introduce the site bond percolation model and prove some auxiliary lemmas concerning this model. The grouping procedure from \cite{KSV}  (see also \cite{KLSV}) is summarized in Section~\ref{grouping}. As in that case, conditioning on a very convenient property of the environment configuration,  where bad ladders are sufficiently spaced, one can develop a multi-scale renormalization procedure. In Section~\ref{sec-renormalization}, we define the renormalized lattices at all levels, as well the notion of open sites and bonds in each level. Section~\ref{prova} is dedicated to the proof of Theorem~\ref{main}. Many details are very similar to the proofs in Theorem of \cite{KSV}, in which case they will be summarized or omitted. We conclude this manuscript in Section~\ref{palavras} with an application in Percolation of Words, a model introduced in \cite{BK} that generalizes ordinary percolation.

\section{Site-bond percolation model on $\mathbb{Z}^2$}
\label{site-bond}

Given the square lattice $G=(\ZZ^2 ,\EE )$ and two parameters $s,p \in
[0,1]$, we define the {\em site-bond percolation model}  as the
Bernoulli percolation model on $\Omega=\{0,1\}^{\VV\cup\EE}$, where sites and bonds are open with probabilities $s$ and $p$, respectively, and independent of each other. Given $u,v\in\ZZ^2$, we use the notation $\{u\leftrightarrow v\}$ to denote the event where there exists a path $\langle u=x_0, x_1,\dots,x_n=v\rangle$ such that the sites $x_i$ and the bonds $\langle x_{i-1},x_i\rangle$ are open for all $i=1,\dots n$ (observe that this doesn't depend if $u$ is open or not); and as usual we define $\{u\leftrightarrow\infty\}=\cup_{v\in\ZZ^2} \{u\leftrightarrow v\}$. Let
$\tilde\theta(s,p)=\tilde P_{s,p}\{0\leftrightarrow\infty\}$ be the
probability of percolation (to avoid confusion we use tilde when referring to the site-bond independent model). Whether $s=p$ we will write only $\tilde P_p$ and $\tilde\theta(p)$. A good reference for the site-bond percolation model is \cite{CS}.

We can now state the following two lemmas.

\begin{lema}\label{bisset} For the Bernoulli site-bond percolation model on $\ZZ^2$,
it holds that the left derivative $\left(\frac{d\tilde\theta(p)}{dp}\right)^-$ is 0 when $p=1$. Consequently, there exists $p^*<1$ such that $\tilde\theta(p)\geq p$ for all $p\geq p^*$.
\end{lema}
\begin{proof}

Let $B_n=\{x\in\ZZ^2;\|x\|_\infty\leq n\}$, $\partial B_n=\{x\in\ZZ^2;\|x\|_\infty= n\}$ and $A_n=\{\omega\in\Omega;0\leftrightarrow \partial B_n\}$.
Since $A_{n} \downarrow\{0\leftrightarrow\infty\}$ %and each event $A_n$ is an increasing event (in the FKG sense). So,
the sequence of functions, $f_n(p)=\frac{\tilde P_p(A_n)}{1-p}$, converges pointwise on $[0,1)$ to the function
$\frac{\tilde\theta(p)}{1-p}$,  as
$n\rightarrow +\infty$. We claim that there exists $p_0<1$ such that this convergence is uniform in the interval $[p_0,1)$. Indeed, if $m\geq n\geq n_0$
\[f_n(p)-f_m(p)\leq \frac{\tilde P_p(A_{n_0}\backslash\{0\leftrightarrow\infty)\}}{1-p}\leq
\frac{1}{1-p}\sum_{i=4n_0}^\infty C^i(1-p)^i.\]
The last inequality was obtained by a Peierls type argument, where $C^i$ is an upper bound for the number of circuits on the dual lattice, formed by $i$ sites or bonds, and containing
the origin (for site-bond percolation on $\ZZ^2$ we can take $C=8$). Therefore, taking $p_0=1-\frac{1}{2C}$, for all $p\in[p_0,1)$ we have that:
\[f_n(p)-f_m(p)\leq\frac{C[C(1-p)]^{4n_0-1}}{1-[C(1-p)]}\leq
C(\frac{1}{2})^{4n_0-2}\]
which shows the uniform convergence of the sequence $(f_n(p))_n$ in the interval $[p_0,1)$, for $p_0=1-\frac{1}{2C}$. It follows that
$\tilde\theta(p)$ is continuous on $[p_0,1)$, since the functions $f_n(p)$ are continuous on $[0,1)$.

Thus, for the left derivative we can write
\begin{align*}
    \frac{d\tilde\theta}{dp}_{\big|_{p=1}}&=\lim_{p\rightarrow
1^-}\frac{\tilde\theta(1)-\tilde\theta(p)}{1-p} \\
&=\ \lim_{p\rightarrow
1^-}\lim_{n\rightarrow
+\infty}\frac{1-\tilde P_p(A_n)}{1-p}=\lim_{n\rightarrow +\infty}
\lim_{p\rightarrow 1^-}\frac{1-\tilde P_p(A_n)}{1-p} \\
&=\ \lim_{n\rightarrow +\infty}\tilde E_1[\#\mbox{ pivotal sites or bonds in the event }A_n]=0.
\end{align*}
The pre-last identity follows from the Russo's Formula and the last one from the fact that for $p=1$, the number of pivotal sites or bonds
for the event $A_n$ is zero, for all $n$. 
\end{proof}

 From the well known fact that $\tilde\theta(\frac{1}{2})=0$ we see that $p^*>\frac{1}{2}$.

 \begin{definicao}
 \label{DR}
 Define $\rho\in(\frac{3}{4},1)$ and for fixed $N>0$, let ${\cal R}_N$ be the set of all rectangles $R$ in $\ZZ^2$ whose sides have at least $N$ sites and at most $6N-1$ sites. Let $D_R$ be the event where there exists an open (site-bond) cluster
contained in  $R$ and satisfying the conditions: {\it i)} it contains on each side of the boundary of $R$ at least $\rho l$ open sites, where $l$ is the number of sites for this side;
 {\it ii)} it contains a circuit made of open sites and open bonds surrounding the center of the rectangle. We define the function $g_N(p)$ as
 \begin{equation}
 \label{g-DR}
g_N(p)=\inf_{R\in{\cal
R}_N}\tilde P_p(D_R).
\end{equation}
\end{definicao}

\begin{lema}\label{cruz1}
For each $\rho\in(\frac{3}{4},1)$, there exists $p(\rho)<1$, such that
$$\lim_{N\rightarrow +\infty} g_N(p)=1$$ for all
$p\geq p(\rho)$.
\end{lema}
\begin{proof} Let $Q_N$ be a square in $\ZZ^2$ made of $N\times N$ sites, $B_{Q_N}^{lr}$ the event that there are at least $\rho N+1$ disjoint paths made of open sites and open bonds, contained in  $Q_N$, and connecting the left side to the right side of $Q_N$. Analogously, let $B_{Q_N}^{tb}$ be the event that there exists at least $\rho N+1$ disjoint paths made of open sites and open bonds, contained in  $Q_N$, and connecting the top to the bottom of $Q_N$. %\marginpar{\tiny ME inside Q}

Equation 2.47 in the Section 2.6 of \cite{Gr} implies that, for Bernoulli bond percolation, there exists $p(\rho)<1$ such that for all $ p\geq p(\rho)$, $$\lim_{N\rightarrow +\infty} \tilde P_{1,p}(B_{Q_N}^a)=1,\ a=lr,tb.$$

With minor modifications, we have the same result (for a different $p(\rho)$) for the site-bond percolation model. Thus, given  $\epsilon>0$ there exists $ N_0< \infty$ such that 

\begin{equation}\label{sobre 7}
\tilde P_p(B_{Q_N}^a)>1-\frac{\epsilon}{7},\ \text{ for all } \ N\geq N_0,\ \text { for all } \ p\geq p(\rho),\ a=lr,tb.
\end{equation}

Given $R\in{\cal R}_{N_0}$, we denote by $l\le L$ the lengths of sides of $R$.  Without lost of generality, we can suppose that the vertices of $R$ have coordinates $(0,0), (0,l),(L,l)$
and $(L,0)$. Let $Q_L$ be the square whose vertices have coordinates $(0,0), (0,L),(L,L)$ and $(L,0)$; for $i=\lfloor\frac{L}{l}\rfloor+1\leq 6$ define $Q_l^j,\ \forall\ j=1,\dots ,i-1$ the squares whose vertices have coordinates
$((j-1)l,0), ((j-1)l,l),(jl,l)$ and $(jl,0)$, respectively, and $Q_l^i$ the square whose vertices have coordinates $(L-l,0), (L-l,l),(L,l)$ and $(L,0)$.

A simple observation shows that if $\frac{3}{4}<\rho$ then $B_{Q_L}^{tb}\cap \left(\cap_{j=1}^i B_{Q_l^j}^{lr}\right)\subset D_R$. Then, by inequality \eqref{sobre 7}, we have that

$$\tilde P_p(D_R)\geq \tilde P_p\left(B_{Q_L}^{tb}\cap \left(\cap_{j=1}^i
B_{Q_l^j}^{lr}\right)\right)\geq 1-\epsilon,\ \forall R\in{\cal R}_{N_0},\ \forall\ p\geq p(\rho).$$
Thus, for all $\epsilon>0$, there exists $N_0=N_0(\epsilon)$ such that for all $N\geq N_0$ it holds $$g_N(p)=\inf_{R\in{\cal
R}_N}\tilde P_p(D_R)>1-\epsilon.$$
\end{proof}

Observe that the event $D_R$ implies the existence of an open cluster in $R$ that contains at least $\rho  \ell$
vertices in each side on $R$, where $\ell$ is the number of sites on that side. Moreover, if this open cluster does not contain the center of $R$, it has a circuit made of sites and bonds surrounding this center. Therefore, if the
center belongs to an infinite open cluster, it must be connected to the net of open crossings of the event $D_R$. This net of open crossings in the event $D_R$ will be called {\em main cluster}.

The next lemma is analogous to Lemma 5.12 in \cite{KSV} and it concerns special crossings in a rectangular region. Let $R$ be the rectangle $[0,L]\times[-\frac{l}{2},\frac{l}{2}]$, $L(R)$ and $l(R)$ be the lengths of the biggest and the shortest side of $R$, respectively. We define the event $C_R=\{(0,0)\overset{ R}{\leftrightarrow}(L,0)\}$ as the set of configurations for which there is a path of open sites and bonds inside $R$, without using bonds on the boundary of $R$, connecting the central vertices of the shortest sides.

\begin{lema}\label{cruz2} There exist $l_0$ and $\varkappa$ positive integers and $0<\bar{p}<1$ close to $1$, such that for all rectangle $R$ with dimensions $l$ and $L$, satisfying $l_0<l\leq L\leq l^2$, it holds that $\tilde{P}_p(C_R)\geq p^{\varkappa }$ for all $p\geq \bar{p}$.
\end{lema} 

\begin{proof}

The proof follows the same steps of  Lemma 5.12 in \cite{KSV}. Define $i=\lfloor\frac{L}{l}\rfloor$ and $r=\lfloor(L-il)/2\rfloor$. Let $A$ be the event where the vertex $(0,0)$ is connected in $[0,l+r]\times [-\frac{l}{2},\frac{l}{2}]$ to at least $\lfloor\frac{l}{6}\rfloor$ vertices in $\{l+r\}\times [-\frac{l}{6},\frac{l}{6}]$ and let $B$ be the event where the vertex $(0,L)$ is connected in $[L-l-r,L]\times [-\frac{l}{2},\frac{l}{2}]$ to at least $\lfloor \frac{l}{6}\rfloor $ vertices in $\{L-l-r\}\times [-\frac{l}{6},\frac{l}{6}]$. We can choose $0<\bar{p}<1$, positive integers $l_0$ and $\varkappa^\prime$ such that $\tilde{P}_p(A)=\tilde{P}_p(B)\geq p^{\varkappa ^\prime},\ \forall p\geq \bar{p}$ and $l\geq l_0$.

For $j=1,\dots,i-2$, let $D_j$ be the event where the vertex $(0,0)$ is connected in $[0,(j+1)l+r]\times [-\frac{l}{2},\frac{l}{2}]$ to at least $\lfloor\frac{l}{6}\rfloor $ vertices in $\{(j+1)l+r\}\times [-\frac{l}{6},\frac{l}{6}]$. Since $\tilde{P}_p(D_j|A\cap D_1\cap\dots\cap D_{j-1})\geq 1-\exp (-c_pl),\ \forall j=1,\dots,i-2$ for some constant $c_p$ bounded away from zero for $p\geq\bar{p}$, we can conclude that (observe that $r+(i-1)l=L-r-l$)
$$\tilde{P}_p(C_R)\geq \tilde{P}_p(A\cap D_1\cap\dots\cap D_{i-2}\cap B)\geq p^{\varkappa }$$
taking $\varkappa >2 \varkappa^\prime$.

\end{proof}

\begin{remark}
\label{rem:cruz2}
Increasing the constant $\varkappa$ if needed, the statement of Lemma~\ref{cruz2} extends uniformly for $p\in[\bar{p},1]$ for any fixed $\bar{p}>p_c$. %\marginpar{\tiny ME corrigi o kappa para varkappa.Entendo que precisamos enfatizar a conexão na direção mais difícil. Mas na aplicação pode ser que o retângulo desejado acabe sendo o oposto. Mas me parece que não precisamos comentar, por ser óbvio. Não?}

%In particular, the statement of the Lemma above holds for all $p>p_g$, since $p_g$ will be chosen latter in such a way that $p_g > p_c$.}
\end{remark}

%\blue{We emphasize that the theorem above holds for all %$p>p_g$, since $p_g$ will be chosen latter in such a way that %$p_g\geq \bar{p}$. } 
We also recall a simple remark that holds for Bernoulli percolation, and which follows at once from coupling. This is Lemma 5.8 in \cite{KSV}. (Of course there is a natural extension if one takes different parameters for site and bond occupation variables, though we don't really need to do this.)

\begin{lema}
\label{sanduiche}
Consider site-bond percolation on a graph $G$. Denote by $P_p$ the probability measure under which all
sites and bonds are independently open  with probability $p$, and let
$\cE$ be some increasing event. If $p_0, p_0^\prime\in [0,1]$ and
$\tilde p = 1-(1-p_0)(1-p^\prime_0)$, then
\begin{equation}
P_p\{\cE\} \ge 1-  (1-P_{p_0}\{\cE\})(1-P_{p^\prime_0}\{\cE\}) \quad \text {for all } p\ge \tilde p.
\label{4.83}
\end{equation}
\end{lema} 
%\subsection{Definition of the parameters}

\section{The environment. Grouping}
\label{grouping}

\medskip
% \marginpar{\tiny ME Bernardo, dá uma olhada em KSV para ver se você concorda com o que está lá e podemos transladar o que for necessário para cá, especialmente para modificar algo na definição 7 aqui. Especialmente lemma 2.3 e lemma 2.4 em KSV}

\noindent Our proof is based on a multi-scale renormalization scheme which depends on a grouping procedure of the bad ladders as introduced in \cite{KSV}. In this section we summarize the main points of this procedure and refer to \cite[Section 2]{KSV} for the full details.

\noindent {\bf Notation.}

\noindent (i) If $C \subset \mathbb{Z}_+$ is a finite set, $\text{span}(C)$ denotes the smallest interval (in $\mathbb Z_+$) that contains $C$;  $\min(C)$ ($\max(C)$) denotes the minimum (maximum, resp.) element of $C$; $\text{diam} (C)=\max\{|x-y|\colon x,y \in C\}$ denotes the diameter of $C$, and $|C|$ denotes its cardinality. (ii) We use $d(A,B)$ to denote the usual Euclidean distance between two sets $A$ and $B$.

\medskip

\noindent {\bf The grouping procedure.} \\
 Given $\xi=(\xi_x)_{x \in \ZZ_+}$, an i.i.d. sequence of Bernoulli random variables with $P(\xi_x=1)=\delta$, the algorithm consists in building a sequence $(\mathbf{C}_k)_k$ of (successively coarser) partitions of $\Gamma=\{x \colon \xi_x=1\}$. To the elements $\cC$ of each partition,  called {\it blocks}, an adequate weight or {\it mass} $m(\cC)$ is attributed and there is also the notion of {\it level} $\ell(\cC)$, which corresponds to the step at which the block is formed (which somehow reflects its complexity). 
 
 The construction depends on a positive parameter $M$ (to be thought as a large integer) and the following are the basic properties, valid for all blocks $\cC$ at each grouping step $k$  %(and also for the limiting partition $\mathbf{C}_\infty$):

\medskip

\noindent (i)\begin{equation}
{\cC}=\text{span} (\cC) \cap \Gamma. \label{2.z}
\end{equation}
\noindent (ii) To each $\cC\in\mathbf{C}_k$,  %we will attribute a {\it mass} $m(\C)=m$ in such a way that
\begin{equation}
d({\cC}, {\cC'}) \ge  L^{\min \{m({\cC}), m({\cC'}),k\}}, \text { for all } 1\leq k \text{ and } \cC, \cC' \in {\mathbf{C}}_k. \label{2.two}
\end{equation}

\noindent (iii) 
% If to each ${\cC} \in {\mathbf{C}}_k $ it holdwe define its {\it level} $\ell({\cC})$, as the smallest $j$ so that $\cC \in {\mathbf{C}}_j$, then
\begin{equation}
\label{2.a}
\ell (\cC) < m(\cC) \text{  for all } \cC \in {\mathbf{C}}_k, \text{ for all } k.
\end{equation}

The construction in \cite{KSV} goes as follows:

\medskip

\noindent {\bf Step 0.} $\mathbf{C}_0=\{\{x\}\colon x \in \Gamma\}$ is simply the partition of $\Gamma$ into
sets of cardinality one. Each of these has mass one. For convenience we label the elements of $\Gamma=\{x_j\}_{j\ge 1}$ in increasing order. 
\medskip

\noindent {\bf Step 1.} At this step one considers all 1-runs of elements of $\Gamma$, i.e. successive elements of $\Gamma$,  $x_i,x_{i+1}.., x_{i+n-1}$ such that $x_{j+1}-x_j< M$, for $j=i,\dots, i+n-2$, and $x_{j+1}-x_j\ge M$ for $j=i-1,j=i+n-1$ (without the first condition in case $i=1$). The blocks of level 1 are formed by these $1$-runs and the mass of such a block coincides with the cardinality of the corresponding run, whose elements were named {\it constituents} in \cite{KSV}. $\mathbf{C}_1$ is formed by the blocks of level 1 and also by those $\{x_j\}$ of $\mathbf{C}_0$ that did not enter any 1-run. \medskip

\noindent {\bf Step $k+1$.}
Having defined the
partitions ${\mathbf C}_{k'}$ for $k' \leq k$ and assuming the validity of
\eqref{2.z}, \eqref{2.two} and \eqref{2.a} when restricted to $k' \leq k$, one considers
 $(k+1)$-runs of blocks in $\mathbf{C}_k$ with mass at least $k+1$: if 
$\cC_1, \cC_2,\dots$ denotes the sequence of all the blocks in $\mathbf{C}_k$ with mass at least $k+1$, labeled in increasing order, the set $r=\{\cC_i,\cC_{i+1},\dots, \cC_{i+n-1}\}$ is said to be a $(k+1)$-run  of length $|r|=n$, where $n \ge 2$, if
\begin{equation*}
d(\cC_{j}, \,  \cC_{j+1})\, < \, L^{k+1}, \; j=i,\dots,i+n-2,
\end{equation*}
and in addition
\begin{equation*}
d(\cC_{j}, \,  \cC_{j+1})\, \ge \, L^{k+1}  \; \;
\begin{cases} \text{for} \; \; j=i-1, j= i+n-1, \; &\text{if}\; i >1 \\
\text{for} \; \; j= i+n-1, \; &\text{if} \; i=1.
\end{cases}
\end{equation*}
A block of level $k+1$ is any set of the form
\begin{equation}
\label{block.k}
\cC= \text{span} (\bigcup_{C \in r} C) \cap \Ga,
\end{equation}
where $r$ is any $(k+1)$-run as above. The mass attributed to $\cC$ is
\begin{equation}
\label{mass.k}
m(\cC)=\sum_{C \in r} m(C) -k(|r|-1).
\end{equation}
The blocks of $\mathbf{C}_k$ that form the $(k+1)$-run in \eqref{block.k} are called {\it constituents} of $\cC$. 

The partition $\mathbf{C}_{k+1}$ is formed by all the blocks of level $k+1$ and all the blocks in $\mathbf{C}_k$ that are not contained in any block of level $k+1$.

\medskip

That the above definition is well set and that properties \eqref{2.z}, \eqref{2.two} and \eqref{2.a} hold at all steps is proven in ~\cite[Section 2]{KSV}, where further properties of the blocks in $\mathbf{C}_k$, valid for all $k$, are also given. We recall only the most important points and refer to \cite{KSV} for a full description and proofs.

The point is that for a suitable relation between $\delta$ and $M$ (see below) the above grouping procedure stops, i.e. with probability one, each $x \in \Ga$ can be incorporated into blocks of higher levels only finitely many times (see Remark \ref{kappa} below). That is, the grouping procedure converges, yielding also a limiting partition that we may call $\mathbf{C}_\infty$ for which \eqref{2.z}, \eqref{2.two} and \eqref{2.a} are also valid for this partition (i.e. can take $k=\infty$ in \eqref{2.two}). It is convenient to demand a bit more, taking the origin as basic reference. This is the content of the following result, which has a crucial role in the renormalization method.

% \begin{lemma}
% \label{lemma111}
% Let us assume that $\de > 0$ and $3 \le L < (64 \de)^{-1/2}$. Under such conditions there exist constants $c_1>0$ and $c_2> \log L$ such that
% \begin{equation}
% \bP \left(\exists\, \cC \in \bigcup_{\l \ge 1}\text{\bf C}_\l
% \colon \text{min}(\cC)=z, m(\cC)=m\right) \le c_1 e^{-c_2 m}, \label{2.5}
% \end{equation}
% for each
%  $m$ and each $z$. In particular we may assume
% \begin{equation}
% c_1(L^m+1)e^{-c_2m} \le 2c_1e^{-c_3m} \label{2.7}
% \end{equation}
% for some constant $c_3 > 0$.
% \end{lemma}
% Before proving the above lemma, we state and prove the following important consequence.

\begin{lema} {\rm(\cite[Lemma 2.4, Lemma 2.8]{KSV})}
\label{chi-variable}
Let $(\xi_x)_{x \in \ZZ_+}$ be an i.i.d. Bernoulli sequence with $P(\xi_x=1)=\delta>0$. Define
\begin{equation}
\chi (\xi) := \inf \left\{ k \ge 0\colon\; \text{min}(\cC) \ge L^{m(\cC)}\; \text{ for all } \cC \in \bigcup_{\ell \ge 1}\mathbf{C}_\ell \text{ with }
 m (\cC) > k \right\},
\label{2.25z}
\end{equation}
with $\chi(\xi)=\infty$ if the above set is empty. 

If  $3 \le M < (64 \delta)^{-1/2}$ we have
\begin{equation}
    \label{chi-eq}
    \PP(\chi < \infty) =1
\end{equation}
and 
\begin{equation}
\label{2.31}
\PP(\chi= 0) > 0.
\end{equation}
\end{lema}

\begin{remark}
\label{kappa-remark}
For each $x \in \Ga$, let $\kappa(x)$ be the random index defined by
\begin{equation}
\label{kappa}
\kappa (x) = \sup\left\{\ell(\cC)\colon x\in \cC \in \bigcup_{0 \leq k <\infty} \mathbf{C}_k \right\}. 
\end{equation}
Since $m(\cC)>\ell(\cC)$ for each $\cC$ as in \eqref{kappa}, an immediate consequence of \eqref{chi-eq} is that $\PP(\kappa(x)< \infty)=1$ under the conditions of Lemma \ref{chi-variable}. 
\end{remark}

\begin{definicao} Given a positive integer $M \ge 3$, the configuration $\xi\in\{0,1\}^{\ZZ_+}$ as above is said to be {\em $M$-spaced} if $\chi(\xi)=0$ for the grouping procedure with scale parameter $M$. 
\end{definicao}

% \begin{defin}
% \label{Lspaced}
% Given $L \ge 3$ an integer, an environment configuration $\gamma$ is said to be \emph{$L$-spaced} if $\chi(\gamma)=0$ for the given choice of scale parameter $L$.
% \end{defin}

% \noindent The construction leading to Lemma \ref{lemma1.1} also yields (see Lemma 2.3 in \cite{KSV})
% \begin{proposition}\label{KSV} Under the conditions of Lemma \ref{lemma1.1} one has
% $$
% \PP(\xi\colon \chi(\xi)<\infty)=1 \text{ and } \PP(\chi(\xi) =0)>0.
% $$
% \end{proposition}

%\begin{definicao}Given a positive integer $M$ and the binary sequences $(\xi_i^a)_{i\in\ZZ},\ a=V,H$, we say that $(\xi_i^a)_{i\in\ZZ}$ is {\em $M$-spaced} if $\xi_i^{a,+}$ and $\xi_i^{a,-}$ are both $M$-spaced semi-infinite binary sequences.
%\end{definicao}
%
%\begin{corolario}\label{delta}Let $\delta <\frac{1}{64M^2}$, then $P_\delta(\xi^V\mbox{and }\xi^H\mbox{are }\mbox{M-spaced})>0$.
%\end{corolario}

\noindent {\bf Applying to our model.}\\
From the independence of $\xi^H$ and $\xi^V$ in our environment configuration $\xi=(\xi^H, \xi^V)$, we have immediately the following 
\begin{corolario}\label{delta}
Under the conditions of Lemma \ref{chi-variable},
\begin{equation}
\label{spaced}
\PP(\xi^{H}, \xi^{V} \mbox{are }\mbox{M-spaced})>0.
\end{equation}
\end{corolario}

%\noindent {\bf Notation:} 
Throughout the rest of the paper we use $\PP_\xi$ to denote the percolation model conditioned on the environment, i.e. $\PP_\xi(A)=\PP(A| \xi)$. 

%\medskip

\begin{remark}
\label{main2}
Notice that $\left\{\xi\colon \PP_{\xi}\{0\leftrightarrow\infty\}>0\right\} $ is a tail event. Thus, we have the following consequence of Corollary \ref{delta}:  
Given $p_g>p_c$ and $p_b>0$, the statement of Theorem~\ref{main} for $\delta$ suitably small follows, if we prove that
we can find $M=M(p_g,p_b)$ finite so that $$\PP_\xi\{0\leftrightarrow\infty\}>0,$$ for any environment configuration $\xi=(\xi^H,\xi^V)$ with $\xi^{H}$ and $\xi^{V}$ $M$-spaced.
\end{remark}

%Theorem  \ref{main} follows once we  prove that for $p_g$ and $p_b$ as in the statement, we can find $M=M(p_g,p_b)$ so that 
%$$\PP_\xi\{0\leftrightarrow\infty\}>0$$
%a.s. on the event $\{\xi^{H} \text{ and } %\xi^{V} \mbox{are }\mbox{M-spaced}\}$.

\section{Renormalization}
\label{sec-renormalization}

%Let $\xi^{H}$ and $\xi^{V}$ be the binary sequences that describe the configuration of good and bad ladders restricted to $\ZZ_+^2$. Given a positive integer $N$ (to be fixed later) we take $\delta$ so as to satisfy the condition of Lemma \ref{lemma1.1} for $M=3N$. Let ${\bf C}_k^H$ and ${\bf C}_k^V$, for $k\in\NN\cup\{\infty\}$, be the corresponding partitions of $\ZZ_+$, obtained by the grouping procedure of the previous section, corresponding to $\xi^H$ and  $\xi^{V}$, respectively. \marginpar{\tiny ME: mencionar ``level $k$}

Given a positive integer $N$ (to be fixed later) and $M=3N$,  let $\xi= (\xi^{H},\xi^{V})$ be $M$-spaced binary sequences, and let ${\bf C}_k^H$ and ${\bf C}_k^V$, for $k\in\NN$, be the corresponding partitions of $\{x \colon \xi^H_x=1\}$ and  $\{x \colon \xi^V_x=1\}$ obtained  at step $k$ of the grouping procedure explained in the previous section.% corresponding to $\xi^H$ and  $\xi^{V}$, respectively.

%\marginpar{\blue{aqui sumiu um parágrafo}}
%Throughout this section we fix a configuration of ladders $\xi= (\xi^{H},\xi^{V})$ with $\xi^{H}$ and $\xi^{V}$ both $M$-spaced. Due to Corollary \ref{main2}, for the proof of the main theorem, it suffices to show that %the conditional probabilities $P_{\xi}(\cdot)=\mathbb{P}(\cdot|\xi)$ satisfy
%\begin{equation}
%\label{percola}
%  \PP_{\xi}\{0\leftrightarrow\infty\}>0,
%\end{equation}
%provided $M=3N$ has been chosen large enough. 

%\medskip

We now describe the renormalization procedure that will be used. It is determined by the partitions $(\mathbf{C}^H_k)_k$ and $(\mathbf{C}^V_k)_k$ and shares some similarities with that used in \cite{KSV}. 

\subsection{The renormalized lattice}\label{renorm}

Given ${\cal C}$ an element of ${\bf C}_k^H$  or  ${\bf C}_k^V$, for $k\in\NN$,  we define its start and end points of ${\cal C}$ as 
$$\alpha({\cal C})=\min\{l\in\ZZ_+;l\in {\cal C}\},$$
$$\omega({\cal C})=\max\{l\in\ZZ_+;l-1\in {\cal C}\}.$$ 

%Analogously, we define $\alpha({\cal C}^H_{\infty,i})$ and $\omega({\cal C}^H_{\infty,i})$ for the clusters of horizontal ladders.

We now define the sites and bonds of our renormalized lattices.

\vskip0.3cm
\noindent \textbf{Step 0:} The {\em sites and bonds at step 0} are simply the sites and bonds of the original lattice $\mathbb{L}^2$.
\vskip0.3cm
\noindent\textbf{Step 1:} Given the partitions ${\bf C}_1^H$  or  ${\bf C}_1^V$, consider each rectangular region between two consecutive blocks of horizontal bonds, say ${\cal C}^H_{1,i}$ and ${\cal C}^H_{1,i+1}$, and two consecutive blocks of vertical bonds,  say ${\cal C}^V_{1,j}$ and ${\cal C}^V_{1,j+1}$, is subdivided into rectangles whose sides have at least $N$ sites (of step 0) and at most $2N-1$ sites (of step 0) in the following manner. Define $l^H_{1,i}$ and $l_{1,j}^V$ as
$$l^H_{1,i}=\left\lfloor\frac{\alpha({\cal C}^H_{1,i+1})-\omega({\cal C}^H_{1,i})}{N}\right\rfloor\geq 3,$$

$$l_{1,j}^V=\left\lfloor\frac{\alpha({\cal C}^V_{1,j+1})-\omega({\cal C}^V_{1,j})}{N}\right\rfloor\geq 3,$$
and for $\ r=1,\dots ,l^H_{1,i}-1$ and $t=1,\dots ,l_{1,j}^V-1$, consider the set of following $l^H_{1,i}\cdot l_{1,j}^V$ rectangles:

$$[(r-1)N+\omega ({\cal C}^H_{1,i}),rN+\omega ({\cal C}^H_{1,i})-1]\times [(t-1)N+\omega ({\cal C}^V_{1,j}),tN+\omega ({\cal C}^V_{1,j})-1],$$
$$[(l^H_{1,i}-1)N+\omega ({\cal C}^H_{1,i}),\alpha ({\cal C}^H_{1,i+1})]\times [(t-1)N+\omega ({\cal C}^V_{1,j}),sN+\omega ({\cal C}^V_{1,j})-1],$$
$$[(r-1)N+\omega ({\cal C}^H_{1,i}),rN+\omega ({\cal C}^H_{1,i})-1]\times [(l_{1,j}^V-1)N+\omega ({\cal C}^V_{1,j}),\alpha ({\cal C}^V_{1,j+1})],$$
$$[(l^H_{1,i}-1)N+\omega ({\cal C}^H_{1,i}),\alpha ({\cal C}^H_{1,i+1})]\times [(l_{1,j}^V-1)N+\omega ({\cal C}^V_{1,j}),\alpha ({\cal C}^V_{1,j+1})].$$

\begin{definicao}[Site of step $1$]
\label{s1}
The {\em site $s$ of step $1$} associated to the rectangle $R$ is the subgraph of $\mathbb{L}^2$, denoted by $\mathbb{L}^2(s)$, formed by all sites and bonds of step $0$ contained $R$. We define $span (s)=R$, the \emph{span of $s$}. Given a site $s$ of step $1$, we define its starting and ending coordinates
\begin{align*}
\alpha^H_1(s)&=\min\{i\in\ZZ_+: \exists j\in\ZZ_+,\ (i,j)\in span(s)\}, \\
\omega^H_1(s)&=\max\{i\in\ZZ_+: \exists j\in\ZZ_+,\ (i,j)\in span(s)\}, \\
\alpha^V_1(s)&=\min\{j\in\ZZ_+: \exists i\in\ZZ_+,\ (i,j)\in span(s)\} \\
\end{align*}
and
\[
\omega^V_1(s)=\max\{j\in\ZZ_+: \exists i\in\ZZ_+,\ (i,j)\in span(s)\}.
\]
\end{definicao}

Observe that each site of step 1 is a rectangular region, made of
sites and bonds of step 0, whose sides have at least $N$ and at
most $2N-1$ sites of step 0, and that all the edges of step 0 in this rectangle belong to good ladders. We consider the natural isomorphism between the set of sites of step 1
and $\ZZ_+^2$. We may associate the origin to the site which contains $(0,0)$ \footnote{Our assumption implies that $\alpha(\mathcal{C}^a_{\infty,1}) \geq 3N$, for $a=V,H$ so that there is one site of step 1 containing the origin.} and then move according to the order in each coordinate. Two sites $s_1$ and $s_2$ of step 1 are said to be adjacent if the previously described natural isomorphism maps them to adjacent sites of
$\ZZ_+^2$.

\begin{definicao}[Bond of step $1$]
\label{bond-s1}
Let $s_1$ and $s_2$ be a pair of adjacent sites of step $1$. There are two possibilities: 
\begin{enumerate}
    \item \emph{(Horizontal bond of step 1)} If $span(s_i)= [c_i,d_i]\times [a,b]$ for $i=1,2$, and $d_1 <c_2$ with $[d_1,c_2]=[\alpha(\mathcal{C}^H), \omega(\mathcal{C}^H)]$ for some $\mathcal {C}^H \in \mathbf{C}^H_1$. We define the horizontal bond of step $1$, $e=\langle s_1,s_2\rangle$, as the subgraph of $\mathbb{L}^2$, denoted by $\mathbb{L}^2(e)$, restricted to the rectangle $[\alpha(\mathcal {C}^H), \omega(\mathcal{C}^H)] \times [a,b]$ and after removal of the bonds
    $$\{\langle (\alpha(\mathcal {C}^H),j);(\alpha(\mathcal {C}^H),j+1)\rangle : j=a,a+1,\dots,b-1\}$$ 
    and
    $$\{\langle (\omega(\mathcal {C}^H),j);(\omega(\mathcal {C}^H),j+1)\rangle : j=a,a+1,\dots,b-1\},$$
    that is, those bonds of level 0 that belong to $s_1$ or $s_2$. 
    
    \item \emph{(Vertical bond of step 1)} If $span(s_i)= [a,b]\times [c_i,d_i]$, $i=1,2$ with $d_1<c_2$, where $[d_1,c_2]=[\alpha(\mathcal{C}^V), \omega(\mathcal{C}^V)]$. The vertical bond of step $1$, $\langle s_1,s_2\rangle$ is defined analogously to the previous case, now considering the rectangle $[a,b] \times [\alpha(\mathcal{C}^V),\omega(\mathcal {C}^V)]$ and removing the bonds of level 0 that also belong to $s_1$ or $s_2$ (See Figure~\ref{fig:elonivel1}).
\end{enumerate} 

\begin{figure}[h!]
%\label{fig:elon1}
\hfil
\includegraphics[width=.9\textwidth]{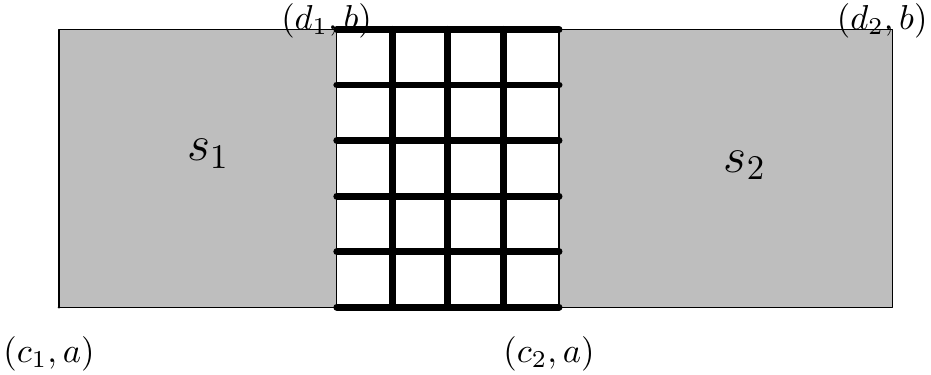}
\caption{The black thick lines correspond to the horizontal bond of step 1.}
\label{fig:elonivel1}
\end{figure}

A bond of step $1$ is said to be \emph{good} if the corresponding block ($\mathcal{C}^H$ or $\mathcal{C}^V$) has mass $1$. Otherwise the bond of step $1$ is said to be \emph{bad}.
%  This horizontal (respectively, vertical) bond of step $1$ is said \emph{good} if the mass of $\mathcal{C}^H$ ($\mathcal{C}^V$) is at most one. If the mass of $\mathcal{C}^H$ ($\mathcal{C}^V$) is larger than 1, we say that this bond of step 1 is \emph{bad}. 
\end{definicao} 

Notice that the subgraph of $\mathbb{L}^2$ that defines any good horizontal bond of step $1$ consists only of a segment of a ladder made of horizontal bonds, whilst for bad horizontal bonds, this subgraph contains horizontal bonds of at least two bad ladders. Similarly for the good vertical bonds of step $1$.

The renormalized lattice made of sites and bonds of step $1$ is isomorphic to $\mathbb{L}^2$ and will be denoted by $\mathbb{L}_{1}^2$. We now define the renormalized lattices of further steps. 

\vskip0.3cm
\noindent \textbf{Step $k (k\geq 2)$:} The sites and bonds of
step $k$ will be defined from the sites and bonds of step $k-1$ as follows:

Given the partitions $\mathbf{C}^H_{k}$ and $\mathbf{C}^V_{k}$, consider each rectangular subgraph of $\mathbb{L}^2_{k-1}$ between two consecutive blocks of horizontal bonds with mass at least $k$, say  ${\cal C}^H_{k,i}$ and ${\cal C}^H_{k,i+1}$, and
two consecutive blocks of vertical bonds with mass at least $k$, 
say ${\cal C}^V_{k,j}$ and ${\cal C}^V_{k,j+1}$. Let us denote by $(s_{m,n})$, with $m=1,\dots,i^*$ and $n=1,\dots,j^*$, the sites of step $k-1$ contained in this subgraph. Since $\xi^{H}$ and $\xi^{V}$ are $3N$-spaced, and each site of step $k-1$ contains at most $(2N-1) \times (2N-1)$ sites of step $k-2$, it holds that $i^*,j^*\geq 3N$, as one can easily verify by induction in $k$.

We subdivide this rectangular subgraph of $\mathbb{L}^2_{k-1}$ at the following manner: define $l^H_{k,i}=\lfloor\frac{i^*}{N}\rfloor,\ l^V_{k,i}=\lfloor\frac{j^*}{N}\rfloor$ and consider the $l^H_i.l_j^V$ rectangles in $\ZZ^2_+$:
\begin{eqnarray} 
\label{retangulo-k}
&&[\alpha^H_{k-1} (s_{(r-1)N+1,1}) , \omega^H_{k-1} (s_{rN,1})]\times [\alpha^V_{k-1} (s_{1,(t-1)N+1}) , \omega^V_{k-1} (s_{1,tN})] ,\\\nonumber
&& [\alpha^H_{k-1} (s_{(l^H_{k,i}-1)N+1,1}) , \omega^H_{k-1} (s_{i^*,1})]\times [\alpha^V_{k-1} (s_{1,(t-1)N+1}) , \omega^V_{k-1} (s_{1,tN})] ,\\\nonumber
&&[\alpha^H_{k-1} (s_{(r-1)N+1,1}) , \omega^H_{k-1} (s_{rN,1})]\times [\alpha^V_{k-1} (s_{1,(l^V_{k,i}-1)N+1}) , \omega^V_{k-1} (s_{1,j^*})] \text{ and }\\\nonumber
&&[\alpha^H_{k-1} (s_{(l^H_{k,i}-1)N+1,1}) , \omega^H_{k-1} (s_{i^*,1})]\times [\alpha^V_{k-1} (s_{1,(l^V_{k,i}-1)N+1}) , \omega^V_{k-1} (s_{1,j^*})], 
\end{eqnarray}
\noindent for $r=1,\dots ,l^H_{k,i}-1\ \ t=1,\dots ,l_{k,j}^V-1$.

For any $(i,j) \in \ZZ^2$, given $R$, any of the rectangles in \eqref{retangulo-k}, we define a {\em site of step $k$} at the following way:

\begin{definicao}[Site of step $k$] The {\em site $s$ of step $k$} associated to $R$ is the subgraph of $\mathbb{L}^2_{k-1}$, denoted by $\mathbb{L}^2_{k-1}(s)$, formed by all sites and bonds of step $k-1$ whose span is contained $R$. We define $span (s)=R$, the \emph{span of $s$}. Given a site $s$ of step $k$, we define its starting and ending coordinates
\begin{align*}
\alpha^H_k(s)&=\min\{i\in\ZZ_+: \exists j\in\ZZ_+,\ (i,j)\in span(s)\}, \\
\omega^H_k(s)&=\max\{i\in\ZZ_+: \exists j\in\ZZ_+,\ (i,j)\in span(s)\}, \\
\alpha^V_k(s)&=\min\{j\in\ZZ_+: \exists i\in\ZZ_+,\ (i,j)\in span(s)\} \\
\end{align*}
and
\[
\omega^V_k(s)=\max\{j\in\ZZ_+: \exists i\in\ZZ_+,\ (i,j)\in span(s)\}.
\]
\end{definicao}

Observe that each site of step $k$ is a rectangular region, made
of sites and bonds of step $k-1$, whose sides have at least $N$ sites of step $k-1$ and at most $2N-1$ sites of step $k-1$. Observe that there is an obvious isomorphism between the set of sites of step $k$ and $\ZZ^2_+$, so we say that two sites of step $k$ are
adjacent, if under this isomorphism they are adjacent sites of
$\ZZ^2_+$. The bonds of step $k$ are defined as in the step $1$ with very minor modifications. 

\begin{definicao}[Bond of step $k$]
\label{bond-sk}
Let $s_1$ and $s_2$ be a pair of adjacent sites of step $k$. There are two possibilities: 
\begin{enumerate}
    \item \emph{(Horizontal bond of step k)} If $span(s_i)= [c_i,d_i]\times [a,b]$ for $i=1,2$, and $d_1 <c_2$ with $[d_1,c_2]=[\alpha(\mathcal{C}^H), \omega(\mathcal{C}^H)]$ for some $\mathcal {C}^H \in \mathbf{C}^H_k$. We define the bond of step $k$, $e=\langle s_1,s_2\rangle$, as the subgraph of $\mathbb{L}^2_{k-1}$, denoted by $\mathbb{L}^2_{k-1}(e)$, restricted to the rectangle $[\alpha(\mathcal {C}^H), \omega(\mathcal{C}^H)] \times [a,b]$ and after removal of the bonds  
    $$\{\langle (\alpha(\mathcal {C}^H),j);(\alpha(\mathcal {C}^H),j+1)\rangle : j=a,a+1,\dots,b-1\}$$ 
    and
    $$\{\langle (\omega(\mathcal {C}^H),j);(\omega(\mathcal {C}^H),j+1)\rangle : j=a,a+1,\dots,b-1\},$$
    that is, excluding those  bonds of level 0 that also belong to $s_1$ or $s_2$. 
    
    \item \emph{(Vertical bond of step k)} If $span(s_i)= [a,b]\times [c_i,d_i]$, $i=1,2$ with $d_1<c_2$, where $[d_1,c_2]=[\alpha(\mathcal{C}^V), \omega(\mathcal{C}^V)]$. The vertical bond of step $k$, $\langle s_1,s_2\rangle$, is defined analogously to the previous case considering the rectangle $[a,b] \times [\alpha(\mathcal{C}^V),\omega(\mathcal {C}^V)]$ and removing the bonds of level 0 that also belong to $s_1$ or $s_2$.
\end{enumerate} 

 A bond of step $k$ is said \emph{good} if the mass of the corresponding block ($\mathcal{C}^H$ or $\mathcal{C}^H$) is at most $k$. Otherwise, the bond of step $k$ is said to be \emph{bad}. 
\end{definicao}

The renormalized lattice made of sites and bonds of step $k-1$ is isomorphic to $\mathbb{L}^2$ and will be denoted by $\mathbb{L}_{k}^2$. Concluding the definition of the renormalized lattices in all scales.

\subsection{Percolation on the renormalized lattice}\label{percrenorm}

\noindent

Since our original model model is a bond percolation model, we define all site of step $0$ as \emph{open} and each bond of step $0$ is \emph{open} is, and only if, it is opened in the original percolation model.

Sites and bonds of further steps will be open or closed depending on certain conditions in the previous steps. Hence, for any step $k\geq 1$, consider the following definitions:

\begin{definicao}[Open sites of step $1$]
\label{s1-open}
A site $s$ of step $1$ is said {\em open} if the event $D$ described in Definition \ref{DR} occurs on the graph $\mathbb{L}^2(s)$; the main cluster (made of sites and bonds of step $0$) defined by event $D$ is called the \emph{skeleton of $s$} and will denote by $skel (s)$. Otherwise the site is said {\em closed}.
\end{definicao}

\begin{definicao} [Open bonds of step $1$]
\label{elo1aberto}
A good bond $e=\langle s_1,s_2 \rangle$ of step 1 is said \emph{open} if the sites $s_1$ and $s_2$ are open and there exists at least one open bond (of step 0) in $\mathbb{L}^2(e)$ connecting the skeletons of $s_1$ and $s_2$,
as described in Definition \ref{DR}; the set of all these open bonds of step $0$ connecting $skel (s_1)$ and $skel (s_2)$ is called the \emph{skeleton of $e$} and will denote by $skel (e)$. Otherwise this (good) bond of step 1 is said {\em closed}.
\end{definicao}

\noindent {\bf Remark.} We shall not try to connect the bad bonds of step $1$, postponing this to a larger scale.
 
\begin{definicao}[Open sites of step $k\geq 2$]
\label{skbom}
A site $s$ of step $k$ is said to be {\em open} if the event $D$  occurs in the subgraph $\mathbb{L}^2_{k-1}(s)$; in this case, we define the \emph{skeleton of $s$} as the union of the skeletons of all sites and bonds of level $k-1$ in the main cluster defined by the event $D$ in $\mathbb{L}^2_{k-1}(s)$. Otherwise this site is said {\em closed}.
\end{definicao}

Observe that for any step the skeletons are made by open sites and bonds of level $0$.

\begin{definicao} [Open bonds of step $k\geq 2$]
\label{elokaberto}
A good bond $e=\langle s_1,s_2 \rangle$ of step $k$ is said \emph{open} if the sites $s_1$ and $s_2$ are open and there exists a path of open sites and bonds of step $0$ connecting the skeletons of $s_1$ and $s_2$; the set of all paths of step $0$ connecting $skel (s_1)$ and $skel (s_2)$ is called the \emph{skeleton of $e$} and will denote by $skel (e)$. Otherwise this (good) bond of step 1 is said {\em closed}.
\end{definicao}

As in the previous steps, we shall not try to connect the bad bonds of step $k$, postponing this to a larger scale. 

\section{Proof of Theorem~\ref{main} }
\label{prova}
As indicated in the Introduction, the proof takes Remark~\ref{main2} into account. Initially, we consider 
the case when $p_g$ is large enough.

\subsection{Restricted crossing. Proof of Theorem~\ref{main} for large $p_g$}
\label{rigorous}

 Lemma \ref{cruz1} implies at once that the probability of a site of step $1$ to be open can be made arbitrarily close to one, provided $N$ will be chosen large enough. Also, if $s_1$ and $s_2$ are two adjacent open sites of step $1$ such that $<s_1,s_2>$ is good, this implies that we find at least $\psi N$ pairs of bad bonds of step 0 with one endpoint (site of step 0) in the skeletons of $s_1$ and the other in the skeleton of $s_2$. Opening at least one of these bonds would immediately connect the skeletons of $s_1$ and $s_2$. Therefore $p_1=1-(1-p_b)^{\psi N}$ is a lower bound for the (conditional) probability of the bond $<s_1,s_2>$ to be open. If in our environment we had only bad layers of mass $1$, a trivial one-step renormalization would bring us to a super critical site-bond model as in Section \ref{site-bond}, with $\tilde p_1=\min\{g_N(p_g), p_1\}$ that tends to one as $N$ tends to infinity. 
 It is also not so difficult to see that if for all $k \ge 1$, every block $\cC$  of mass $k$  would be a run of $k$ consecutive ones in the $\xi$ sequence, we would apply the same type of reasoning at all scales or steps and would conclude that there is percolation for such environments. This includes a class of  hierarchical type of environments.  Developing the previous argument at all scales (or steps), and considering the behavior of  $g_N(p_k)$ where $p_k$ quickly tends to 1, e.g. as $p_k=1 - (1-p_{k-1})^{\psi N}$ the statement is obtained. Full details in this case are given in Section 2.2.2 of \cite{Lima}.
%  \red{(Figuras são bem-vindas)}
  
 Nevertheless,  our random environment $\xi$ presents us with a substantially more involved situation. There are basically two main differences:
 \begin{itemize}
     \item {} Even for simple blocks $\mathcal{C}$ of level 1, for which its mass $m(\mathcal C)=k$ coincides with its cardinality, we should observe that the bad ladders are not all consecutive. This is not a true problem, but we need to connect a pair collinear points in the corresponding skeletons of $s_1$ and $s_2$, fixed on opposite sides of the good bond $<s_1,s_2>$ of step $k$ with a controllable probability cost. For that we proceed by using only order $\sqrt{N}$ collinear pairs which are well separated so that Lemma \ref{cruz2} can guarantee the connections using disjoint sets of edges and within the graph that defines the edge. Proceeding in this manner we can again reproduce a small variation of the argument used for the previous case, at the expense of a factor of order $1- (1- {p_g}^{\kappa (k-1)})^{ N^{k/2}}$ for each pair $s_1$ and $s_2$ of open adjacent sites of step $k$. Since this grows to $1$ fast enough, it is possible to extend the previous argument without too much pain. 
     
     \item{} Far more complicated is the fact that, depending on its level, the cardinality of block $\mathcal{C}$ can indeed be much larger than its mass. The motivation behind that definition is that a larger good interval in the middle allows for better growth, and this explain the mass does not count all the layers that correspond to an $i \in \mathcal{C}$, i.e. for which $\xi_i=1$. This forces us to find a much better description of the crossings through the corresponding bond that corresponds to $\mathcal {C}$. On the other hand, recall that in our graph, we have removed each full line of orthogonal edges in the middle of two consecutive bad ladders. This not only makes our main result stronger, but justifies the way we handle the difficulty. When passing through a bad ladder of level 1, a path needs to move as described in the previous paragraph i.e. crossing consecutive bad bonds in a straight manner. We use this at all steps, giving space between the paths under consideration when crossing good layers of any step. Thus, when discussing percolation through a bond (say horizontal) which corresponds block $\mathcal {C}^H$ of mass $m$ and some level $\ell <m$, we may indeed use only ``directed" movements at the proper scales. This brings us to the context treated in \cite{KSV}, and we may indeed use the recursion developed there, which involves a suitable decomposition of $\mathcal{C}^H$. The construction  developed in sections 6 and 7 of that paper, and which consists of a very specific procedure for checking connections, works here as well. This will be described next.
 \end{itemize}

We now establish conditions on the scale parameter $N$. Let $\rho$, $\psi$ and $p^\prime$ be such that
$\frac{3}{4}<\rho<1$, $\psi=2\rho-1>\frac{1}{2}$ and $p^\prime:=\max\{p^*,p(\rho)\}$, where $p^*$ and $\ p(\rho)$ are given by Lemmas~\ref{bisset} and~\ref{cruz1}, respectively.

We can now define $\p_0$ (the lower bound for $p_g$ in this proof) as $\p_0=\max\{p^\prime, p^{\prime\prime}\}$, where $p^{\prime\prime}$ will be defined in Lemma \ref{estim} below and $\p_k=1- (1-\p_{0})^{k+1}$ for all $k \ge 1$. Given such $\p_0<1$, Lemma \ref{cruz1}, gives $N^\prime$ large enough so that
\begin{equation}
    \label{03mar-a}
g_N(\p_0)>1-(1-\p_0)^2,\ \forall\ N\geq N^\prime
\end{equation}
and $N^\prime\geq\frac{l_0}{9}$, where $g_N$ is the function defined Definition~\ref{DR} and $l_0$ is given by Lemma \ref{cruz2}. Using Lemma \ref{sanduiche} we have at once that for all $N \geq N'$
\begin{equation}\label{inducao}
g_N(\p_k)\geq \p_{k+1},\forall k\geq 0.
\end{equation}

%The parameter $\widehat{\delta}(p_b,p_g)=\frac{1}{64{(3N)}^2}$, $\forall\ p_b>0,\
%\forall\ p_g>\widehat{p}$, observe that $N$ is a function of $p^\prime$ and $p_b$.

%For all $\delta<\widehat{\delta}(p_b,p_g)$, Corollary \ref{delta} holds that $\xi^V\mbox{and }\xi^H\mbox{are }\mbox{3N-spaced}$ with positive $P_{p_b,p_g,\delta}$-probability. Then, as stated in Corollary \ref{main2}, it is enough to show that there exists $c>0$ such that
%\begin{equation}\label{c}
%P_{p_b,p_g,\delta}(0\leftrightarrow\infty| (\xi^i_V)_{i\in\ZZ},(\xi^i_H)_{i\in\ZZ}))>c>0,
%\end{equation}
%for all configuration of ladders $(\xi^i_V)_{i\in\ZZ},(\xi^i_H)_{i\in\ZZ})$ $3N$-spaced. 

%Defined the parameters, we perform a multi-scale renormalization, where sites and bonds of level $k$ will be open as defined in Section~\ref{percrenorm}. 

\begin{lema} {\bf (Main estimate)} 
\label{main-estimate}
Let $p_b>0$ and $p_g=\p_0$ as defined above. There exists $M_0=M_0(p_g, p_b)$ such that for all $M\geq M_0$ whenever $\xi=(\xi^H, \xi^V)$ is an $M$-spaced environment configuration we have: 
\noindent (i) For any $k \ge 0$ and any site $s$ of step $k$
\begin{equation}
\label{lemma-main-s}
\PP_{\xi}\{ s \text{ is open }\}\ge \mathfrak{p}_k, 
\end{equation}
where $\p_0=p_g$ and for all $k \ge 1$, $\p_k=1- (1-\p_{0})^{k+1}$.

\noindent (ii) If $s_1$ and $s_2$ are two adjacent sites of step $k$ and $<s_1,s_2>$ is a good bond of step $k$, then  
\begin{equation}
    \label{lemma-main-b}
    \PP_{\xi}\{(<s_1,s_2> \text{ is open }| s_1 \text{ and } s_2 \text{ are open })\} \ge {\p}_k.
\end{equation}
\end{lema}
 %\marginpar{\tiny ME: depois vemos se usar $q_k$ ou algo como $p_k$. Isso porque como iremos dar uma referência à prova em ksv e lá $q_k=1-p_k$ fica meio enrolado.} 
\noindent {\bf{Proof.}} 
We choose $M_0=3N$, where $N$ is given by equation~\ref{24-03-22-1} below (indeed, $N$ is a function of $p_b$ and $p_g$). There is nothing to prove when $k=0$. For sake of discussing the general argument, we start by considering the step $k=1$, which will give us a lower bound for $M_0$, and then we will proceed with the induction step. 
 Indeed, according to Lemma \ref{cruz1}, we see that if $s_1$ is a site of step $1$,
then (recall that $M\geq 3N$)   
 \begin{equation}
     \PP_{\xi}\{ s \text{ is open }\}\ge g_N(p_g) \ge 1- (1-p_g)^2, 
 \end{equation}
 provided $M$ is taken as indicated above (see \eqref{03mar-a}). It is also simple to see that \eqref{lemma-main-b} holds for $k=1$. Indeed at this step, while conditioning on the two sites of step 1 being open, we have at least $\psi N$ edges $<x_1,x_2>$ of step 0, and where $x_i \in {skel}(s_i)$ for $i=1,2$. It suffices to have one of these edges to be open. Therefore the left hand side of \eqref{lemma-main-b} is bounded from below by $1- (1-p_b)^{\psi N}$ and, as already mentioned, it suffices to take $N$ so that 
 $(1-p_b)^{\psi N} \leq (1-p_g)^2$. %\marginpar{\tiny ME ajustar com o que foi dito na secao da definição dos parâmetros. Colocar aqui a escolha de $M$.}

To achieve the induction step, notice that the validity of \eqref{lemma-main-s} and \eqref{lemma-main-b} at step $k$  implies that of \eqref{lemma-main-s} at step $k+1$ as one sees from \eqref{inducao}. Thus the main point in the induction step consists in the proof of \eqref{lemma-main-b} at step $k+1$.

At this point we follow the argument given in \cite{KSV} since our problem can be translated to that situation, as we now explain. To fix notation, and without loss of generality we consider $<s_1,s_2>$ as a good horizontal bond of step $k$, i.e. it corresponds to a block $\cC=\cC^H$ with $m(\cC)=k$, and we think of $s_1$ located to the left of $s_2$.
 
In \eqref{lemma-main-b} we are conditioning on the two adjacent sites $s_1$ and $s_2$ being open, and speak about the probability to have $skel(s_1)$ and $skel(s_2)$ connected by an open path through $<s_1,s_2>$. Translating it to the   construction in \cite{KSV}, this corresponds to 
\begin{equation}
\label{drilling-1}\PP_\xi( {\Psi^k}  \leftrightsquigarrow_{_{\! \! \! \! \! \!\! \!
\! \! \! \!\!\!  B(\cC)}} {\Upsilon^k})\ge \p_k, 
\end{equation}
where $\Psi^k$ represents the set of sites (of step 0) on the right vertical side of $skel(s_1)$ and analogously $\Upsilon^k$ represents the set of sites (of step 0) on the left vertical side of $skel(s_2)$. They are therefore separated by the bond $<s_1,s_2>$ which is defined on a segment of the vertical layer defined by a block $\cC$ with $m(\cC)=k$, and which we denote as $B(\cC)$. We recall that only the portion of $B(\cC)$ contained in $<s_1,s_2>$ is allowed in  \eqref{drilling-1}. The precise definition of this event is given below; it follows \cite{KSV} and corresponds to the notion of $\Psi^k$ and $\Upsilon^k$ being {\em chained} introduced therein.

As already mentioned, the main difficulty is that a layer that corresponds to a block of mass $m$ can include a much larger number of single bad layers, depending on its level. To deal with this, and to provide an induction estimate, in \cite{KSV} one first makes a more detailed description of a block $\cC$ of a given mass $m$ and level $\ell$. A preliminary fact is that a block $\cC$ of mass $m$ and level $\ell$ has at most $m-\ell+1$ constituents, as one easily sees. Also, the interval between two constituents of an $\ell$-run, also called {\em porous medium of level $\ell$} contains only blocks of smaller masses.  \cite[ Lemma 6.4]{KSV} says that any block $\cC$ of mass $m\ge 2$ has a representation, thereby called {\it descending decomposition}, described as follows: if $\cC$ has mass $m$ one can always find a sequence of integers 
 $\min(\cC)=f_1<g_1 <f_2<g_2<...<f_v<g_v \leq \max(\cC)-1$, in such a way that if one restricts the environment to the interval $[f_s,g_s]$, then $\Gamma \cap [f_s,g_s]$ is a unique block, call it $\tilde \cC_s$, and  
 $m (\tilde \cC_s)=\tilde m_s$, where $m-1=\tilde m_1>\tilde m_2>...>\tilde m_v$ and the intervals $[g_{s-1}, f_s-1]$ are porous media of level $\tilde m_s$ for the given environment configuration, with 
 $$
M^{\tilde m_s} \leq f_s-g_{s-1} \leq M^{\tilde m_s +1}; \\
 \max(\cC)-M <g_v, [g_v+1, \max(\cC)-1] \cap \Gamma=\emptyset
 $$
 This collection $\tilde m_1,... \tilde m_v$, also called \emph{itinerary}, provides a proper way to proceed inductively so as to build an open path crossing the corresponding bad layers, as done in \cite[Section 7]{KSV}.  % telling us which scales we have to examine wile trying to connect sites separated by a bond of step $m$ that corresponds to such a $\cC$. 

The proof (by induction) of the estimate \eqref{drilling-1} can be obtained here exactly in the same manner as that of \cite[Theorem 5.15]{KSV}. It corresponds to the property thereby denoted by $(b_k)$ and is obtained with the help of a more detailed description through properties $(b_k)'$ and $(b_k)''$  in \cite[Proposition 7.1]{KSV}. Our model is not oriented, but since we have removed the edges that are perpendicular to two consecutive bad layers, it is not so surprising to use the same type of restricted paths: the paths move straight within bad layers, in good regions one uses Lemma \ref{cruz2} with care to consider disjoint regions so as to keep the needed independence used for the probability estimates. In order to achieve these estimates, one looks at suitable crossings at proper scales, and consider suitable subsets of $\Psi^k$ and $\Upsilon^k$, forming hierarchical collections $\overline{\Psi^k}=\{\overline{\Psi^j}\}_{0 \leq j \leq k}$ and $\overline{\Upsilon^k}=\{\overline{\Upsilon^j}\}_{0 \leq j \leq k}$ of subsets of $\Psi^k$ and $\Upsilon^k$, respectively. For each $j$, these subsets are formed by  elements of $\Psi^k$ ($\Upsilon^k$) that belong to the sites of scale $j$ ($0 \leq j \leq k$) contained in $skel(s_1)$ ($skel(s_2)$, respectively). The important thing is that, given that $s_1$ and $s_2$ are open, each time we go down the scale (from $j$ to $j-1$),  we can select only pairs that are well separated in the vertical coordinate, (at distance of order $\sqrt{N}$) and will have enough such pairs, at least $J:=\lfloor \frac{\psi}{3} \sqrt N\rfloor$ of them. This set will be denoted by $\mathcal{M}(\Psi^k)$ and $\mathcal{M}(\Upsilon^k)$.
% which form {\em matching pairs} at all scales, in the sense that they correspond to the same vertical coordinates (in this case of an horizontal edge).

The following corresponds essentially Definition 6.13 in \cite{KSV}.

\begin{definicao} \label{chained}

When $\cC$ has mass $m$, level $\ell$ and $k \in\{\ell-1, ..., m-1\}$ the definition for $\overline{\Psi^k}$ and $\overline{\Psi^k}$ being chained through the horizontal edge determined by $\cC$ means:

\noindent (i) When $k=0$, we must have $\ell=1$ and $m$ is exactly the total number of bad layers of step 0 contained in $\cC$ we might separate in two cases:

$\bullet$ If $\cC$ is a sequence of $m$ bad layers of step 0. Recall that $\Psi^0$ and $\Upsilon^0$ are sites of step 0 which have the same vertical coordinates. The property of $\Psi^0$ and $\Upsilon ^0$ being chained corresponds simply to the existence of a straight (horizontal) open path connecting them.  

$\bullet$ If $\cC$ has a certain number of constituents (each of them being as in the previous case). We then demand the existence of straight horizontal paths within each of the bonds determined by the constituents of $\cC$ and that in between them (good layers) we have an open path restricted to a strip of width $\sqrt{N}$. %\footnote{For the probability estimate we apply here Lemma \ref{cruz2}.} 

\noindent (ii) When $k \ge 1$, we still have to distinguish the case $k \ge \ell$ from that when $k = \ell - 1$:  (a) In the first case, we have to consider $J$ pairs $\overline{\Psi^{k-1}}$ and $\overline{\Upsilon^{k-1}}$ of separated matching pairs of scale $k-1$ and simply require that the property of being chained is valid for one of them. (b) When $k=\ell-1$ we need to consider the way $\cC$ is formed, i.e. 
considering its constituents $\cC_1, \dots, \cC_r$, their masses $m_v$ and levels $\ell_v$, with $v=1,...r$. We then require the existence of open sites of step $k$
$s'_1,\dots, s'_r$ such that: (b1) each $s'_v$ and $s'_{v+1}$ form a matching pair with respect to the horizontal edge determined by $\cC_v$; (b2) $\overline{\Psi^k}$ and $\overline{\Upsilon^k}(s'_1)$ are chained through the edge determined by $\cC_1$, $\overline{\Psi^k}(s'_v)$ and $\overline{\Upsilon^k}(s'_{v+1})$ are chained through the edge determined by $\cC_v$, for $v=2, ..., r-1$, and $\overline{\Psi^k}(s'_r)$ and $\overline{\Upsilon^k}$ are chained through the edge determined by $\cC_r$; (b3) For each $v=1, \dots, r$, there is an open path of step $k$, from $s'_v$ to $s'_{v+1}$. This path is contained in a strip of order $\sqrt N$
in terms of the scale for step $k$ i.e. of order $N^{k+\tfrac{1}{2}}$.
\end{definicao}

 \noindent {\bf Notation.} Let $\varkappa>0$ given by Lemma~\ref{cruz2} and for $m\ge 1$ 
 we recursively set:
\begin{eqnarray}
\label{6N.1}
\nonumber p_{0,m}&:=&p_{b}^m p_{g}^{\varkappa(m-1)},\\
p_{j,m}&:=&(1-(1-p_{j-1,m})^J)\p_j^{\varkappa(m-j-1)},\qquad 1\le
j\le m-1,\\\nonumber p_{m,m}&:=&1-(1-p_{m-1,m})^J,\\\nonumber
\end{eqnarray}
where $J=\lfloor{ \frac\psi3N^{1/2}}\rfloor$, with $\psi=2\rho -1$ as introduced before, and we recall
 $\p_j=1-\q_j$, $\q_j=\q_0^{j+1}$, $\p_0=p_{g}, \q_0=1-p_g$, as in   Theorem \ref{main}

Taking the previous definition into consideration, the same proof used in \cite[Proposition 7.1]{KSV} gives us that the following two statements hold true at all steps $m$, where the notion of hierarchical sets refer to the one given above, assuming two adjacent open sites of step $m$, $s_1$ and $s_2$, with the edge $<s_1,s_2>$ that corresponds to a layer $\cC$, of mass $m$. The proof is probably one of the most delicate parts of \cite{KSV}. The detailed structure of a block $\cC$ as its descending decomposition described above plays an important role. Another important ingredient is the way to reconstruct the probabilities $p_{j, m+1}$ from those of scale $m$, through the algorithm described as {\em checking procedure} in \cite[Definition 7.4]{KSV}.

\begin{proposition}
\label{chained-estimates}
Under the conditions of Theorem \ref{main}. and the just described setup, i.e. given two adjacent open sites of step $m$, $s_1$ and $s_2$, with the edge $<s_1,s_2>$ that corresponds to a layer $\cC$, of mass $m$, we  have the validity of $(b_m)^\prime$ and $(b_m)^{\prime\prime}$ as defined below. 

\noindent For $m \ge 1$: 

\noindent $(b_m)^\prime$  For every block $\cC\in \mathbf{C}_\ell$ of  mass $m$ and level $\ell$,  every  $j \in\{\ell-1,\dots, m-1\}$  and every choice of
hierarchical $j$-sets $\overline {\Psi^j},\overline { \Upsilon^j}$ defined as above,  one has
\begin{equation}
\PP_\xi (\overline {\Psi^j}  \leftrightsquigarrow_{_{\! \! \! \! \! \!\! \!
\! \! \! \!\!\!  B(\cC)}} \overline {\Upsilon^j})\ge p_{j,m}.
\label{6N.2}
\end{equation}

\noindent For $m \ge 2$:

\noindent $(b_m)^{\prime\prime}$ For every $B(\cC)$, $j$ and %\marginpar{\tiny ME ainda para decidir o quanto vamos colocar do que consta como def. 6.14 em KSV e que entraria logo acima disto aqui. A notação $\mathcal{M}(\Upsilon)$} 
$\overline {\Psi^j},\overline { \Upsilon^j}$ as in $(b_m)^\prime$,
let $\mathcal{M}(\overline{\Psi^j})$ (and similarly $\mathcal{M}(\overline{\Upsilon^j})$)  be the collection of the $J$ checked subsets of $\overline{\Psi^j}$ in each scale; and every $s\in\{0,\dots,j-1\}$, the conditional distribution (under $\PP_\xi$)  of the number of
$\Upsilon^s_{\langle {\mu_{\langle j,s+1 \rangle}},i \rangle}\in
\mathcal{M}(\overline{\Upsilon ^j})$ that are chained to $\overline
{\Psi^j}$, given that $\Upsilon^{s+1}_{{\mu_{\langle j,
s+1\rangle }}}$
 is chained to $\overline {\Psi^j}$, where $F_p$ denotes
   the distribution of a Binomial random variable with $J$ trials and success probability $p$,
    conditioned to have at least one success. That is,
\begin{eqnarray*}
 | \{i \colon
\Upsilon^s_{\langle {\mu_{\langle j,s+1 \rangle}},i \rangle}\in
\mathcal{M}(\overline{\Upsilon ^j}) \colon \Upsilon^s_{\langle
{\mu_{\langle j,s+1 \rangle}},i \rangle}
 {\text {chained to }} \overline {\Psi^j} \}| \Big|
\big[   \Upsilon^{s+1}_{{\mu_{\langle j, s+1\rangle }}}  {\text {
chained to }} \overline {\Psi^j} \big] \\\succeq
F_{p_{s,m}},\label{6N.3}
\end{eqnarray*} 
with $\succeq$ standing for stochastically larger in the usual
sense.
\end{proposition}

We refer to \cite{KSV} for the proof of this proposition.
\bigskip

\noindent {\bf Conclusion of the induction step in \eqref{lemma-main-b}.}

Back to Lemma \ref{main-estimate} we need to see how to proceed the induction step in \eqref{lemma-main-b} from the estimates in Proposition \ref{chained-estimates}. This involves the same  checking procedure as in  the proof of the previous proposition, together with the choice of our initial parameters. The details are very similar to the corresponding result in \cite{KSV}.

Taking into account the previous arguments and the notation introduced in \eqref{6N.1}, it suffices to check that for all $m \ge 2$, 
\begin{equation}
    \label{24-03-22-1} 
    C(1-p_{m-1, m})^{\frac{\psi}{6} N} \leq \q_{m},
\end{equation}
where $C$ is a fixed universal constant.

Indeed, and since $N$ has been taken large enough, it suffices to prove that 
\begin{equation}
    \label{24-03-22-2}
p_{m,m} \geq \p_{m} 
\end{equation}
for all $m \ge 2$.  
To check this, let 
$$
\Theta= \prod_{k=0}^\infty \p_k >0
$$
which is an increasing function of $p_g$.  

One goes back to the checking procedure described in the previous proof, where now we have replaced $m+1$ by $m$ for notational convenience. Now we leave the testing with probability $\tilde p:=p_b^m$ for the final step of the $0$-boxes. Combine this with the trivial observation that if one has $t$ independent Bernoulli random variables with success probability $p\tilde p$, then the probability of no success is bounded from above by
$$
(1 - \tilde p)^{\lfloor {t}p/2\rfloor} + e^{-{t}I_p(p/2)}
$$
where $I_p(x):=x \log(x/p) +(1-x) \log(1-x)/(1-p)$ for $x \in (0,1)$ and $t$ is a positive integer. 
Following the previous checking procedure at all steps $i$ from 0 to $m-2$, each $i$-box is tested (independently of anything else) with probability $\p_i^{\varkappa(m-i-1)}$, and at the end the $0$-box has to be approved with probability $\tilde p$, we end up with
\begin{align}
\label{24-03-22-3a}
1-p_{m,m} &\leq (1-p_b^m)^{ 4(J/2)^m \prod_{i=0}^{m-2} \p_i^{\varkappa(m-i-1)}} \\
\label{24-03-22-3}
&+ \sum_{i=1}^{m-1} \exp[{-4(J/2)^{i+1}\prod_{j=2}^i \p_{m-j}^{\varkappa(j-1)} f(\p_{m-i-1}^{i\varkappa})]},
\end{align}
where 
$$
f(p)=I_p(p/2)=\left(1-\frac{p}{2}\right)\log \left(\frac{2-p}{1-p}\right)- \log2.
$$
One sees that  $N_0=N_0(p_g, p_b)$ large can be taken so that for all $m \ge 2$ and all $N \ge N_0$, the term \eqref{24-03-22-3a} is bounded from above by 
$\frac12\q_m$. 
The lemma below controls the term \eqref{24-03-22-3}.

\begin{lema}\label{estim}For the sequence $(\p_k)_{k\in\NN}$ defined above, it holds the inequality 
\begin{equation}\label{gd}
\sum_{i=1}^{m-1}
\exp\left[-4\big(\frac{J}{2}\big)^{i+1}\prod_{j=1}^{i-1}\p^{j\varkappa}_
{m-j-1}f(\p^{i\varkappa}_{m-i-1})\right]
\leq \frac{1-\p_m}{2},\ \forall k\in\NN.
\end{equation}
\end{lema}
\begin{proof}
This proof follows the same steps of the end of Section 7 in \cite{KSV}. We have that

\begin{align*}
&\sum_{i=1}^{m-1}
\exp\left[-4\big(\frac{J}{2}\big)^{i+1}\prod_{j=1}^{i-1}\p^{j\varkappa}_
{m-j-1}f(\p^{i\varkappa}_{m-i-1})\right] \\
%\leq\sum_{i=1}^{k-1}\exp \left[-C_N^{i+1}\left(\prod_{j=1}^{i-1}\frac{r^{j\kappa}_{k-j-1}}{2}\right)f(r^{i\kappa}_{k-i-1})\right] \\
%\label{quebraem2}
&=\ \sum_{i=1}^{\frac{m}{2}-1}\exp \left[-4\big(\frac{J}{2}\big)^{i+1}\prod_{j=1}^{i-1}\p^{j\varkappa}_{m-j-1}f(\p^{i\varkappa}_{m-i-1})\right]\\
&\ + \sum_{i=\frac{m}{2}}^{m-1}\exp \left[-4\big(\frac{J}{2}\big)^{i+1}\prod_{j=1}^{i-1}\p^{j\varkappa}_{m-j-1}f(\p^{i\varkappa}_{m-i-1})\right]
\end{align*}

The second term of the r.h.s. can be bounded above by
\begin{equation}\label{igrande}
\frac{m}{2}\exp \left[-4\left(\frac{J}{2}\right)^{\frac{m}{2}}(\prod_{j=1}^\infty \p_j)^{m\varkappa}f(\p_0^{(m-1)\varkappa})\right]\leq \frac{1-\p_m}{4}
\end{equation}
since our choice of $p^{\prime\prime}$ close to one; for the first term, we can give the upper bound

\begin{equation*}
\frac{m}{2}\exp \left[-2J(\prod_{j=1}^\infty \p_j)^{\varkappa}f(\p_{\frac{m}{2}}^{\frac{m}{2}\varkappa})\right]\leq
\frac{m}{2}\exp \left[-4f(\p_{\frac{m}{2}}^{\frac{m}{2}\varkappa})\right]
\end{equation*}
in the last inequality we are also taking $p^{\prime\prime}$ close to one and assuming that $J>1$ (that is, $N>36$).

Then, it is enough to show that
\begin{equation}
\frac{m}{2}\exp \left[-4f(\p_{\frac{m}{2}}^{\frac{m}{2}\varkappa})\right]\leq\frac{1-\p_m}{4}
\end{equation}
or equivalently
\begin{equation}
4f(\p_{\frac{m}{2}}^{\frac{m}{2}\varkappa})\geq -(m+1)\ln (1-p^{\prime\prime}) + \ln (2m)
\end{equation}
using the definition of $f$, we need to show that
\begin{equation}\label{1estim}
-\ln (1-\p_{\frac{m}{2}}^{\frac{m}{2}\varkappa})\geq\frac{1}{4}\left[\ln (4m)-(m+1)\ln (1-p^{\prime\prime})\right]
\end{equation}

Observe that

\begin{align}
1-\p_{\frac{m}{2}}^{\frac{m}{2}\varkappa} &=1-\sum_{i=0}^{\frac{m\varkappa}{2}}\binom{\frac{m\varkappa}{2}}{i}(-1)^i\exp\left[i(\frac{m}{2}+1)\ln (1-p^{\prime\prime})\right] \\
&=\ \sum_{i=1}^{\frac{m\varkappa}{2}}\binom{\frac{m\varkappa}{2}}{i}(-1)^{i+1} (1-p^{\prime\prime})^{i(\frac{m}{2}+1)} \\
&\leq\ \sum_{i=1}^{\frac{m\varkappa}{2}}\binom{\frac{m\varkappa}{2}}{i} (1-p^{\prime\prime})^{i(\frac{m}{2}+1)}\leq 2^{\frac{m\varkappa}{2}} (1-p^{\prime\prime})^{(\frac{m}{2}+1)}
\end{align}
Taking logarithms into the equation above

\begin{equation}
-\ln (1-\p_{\frac{m}{2}}^{\frac{m}{2}\varkappa})\geq -\frac{m\varkappa}{2}\ln 2-(\frac{m}{2}+1)\ln (1-p^{\prime\prime})\geq -\frac{1}{4}(m+1)\ln (1-p^{\prime\prime})
\end{equation}
In the last inequality, we are using again that $p^{\prime\prime}$ close to one to show that $\frac{\varkappa}{2}\ln 2<-\ln (1-p^{\prime\prime})$. This proves equation \ref{1estim} and finishes the proof of Lemma \ref{estim}.
\end{proof}

Recalling Remark~\ref{main2}, this gives the proof of Theorem~\ref{main} for $p_g$ large enough and $p_b>0$.

\subsection{Proof of Theorem~\ref{main} for  $p_g>p_c$} 

Recalling Remark~\ref{main2}, it suffices to prove:
%he extension for $p_g>p_c$ follows essentially the same steps contained in Section 8 of \cite{KSV}. We would like to prove the following statement:

\begin{teorema}
 For any $p_b>0$ and $p_g>p_c$, there exists $M_0(p_g,p_b)$ finite so that for all $M\geq M_0$, $$\PP_\xi\{0\leftrightarrow\infty\}>0,$$ for any environment configuration $\xi=(\xi^H,\xi^V)$ with $\xi^{H}$ and $\xi^{V}$ $M$-spaced.
\end{teorema}
\begin{proof}
[Outline of the proof]
The extension for $p_g>p_c$ follows essentially the same steps contained in Section 8 of \cite{KSV}. It is enough to extend Lemma~\ref{main-estimate} for any $\p_0>p_c$. The central point is to modify the renormalized lattice at step 1 in such a way that each site and bond of step 1 will be open with probability at least $\hat{p}=\max\{p^\prime, p^{\prime\prime}\}$; therefore, the same argument as the proof of Lemma~\ref{main-estimate} follows.

We summarize the main points of this modification, full details can be found in Section 8 of \cite{KSV}.

Replace $\rho$ by $\hat{\rho}\in (0,\theta(p_g))$.
We modify the definitions of site and bond of step $1$. Constructing the renormalized lattice of step $1$ as described in Section~\ref{renorm}, the graph that defines a site of step 1 is taken as the subgraph contained in a centralized square of side $0{.}9N$, and the graph that defines the bond of step $1$ will incorporate the external layer of width at least $0{.}1N$ between the new sites of step $1$ that defines its end-vertices.

Each site of step $1$ is open if its skeleton has at least $0{.}9\hat{\rho}N$ sites of step $0$ on each side of its boundary; each bond of step $1$ is open if there exists an open path of bonds of step $0$ connecting the skeletons of its end-sites of step $1$; observe that all sites and bonds of step $1$ are open independently. Increasing $N$ if needed, we see that each site and bond of step is open with probability at least $\hat{p}$. Finally, take $M_0=3N$.

\end{proof}

\section{An application}
\label{palavras}

% As explained in the previous sections, the main result of this paper uses the techniques developed in \cite{KSV} in combination with those that were introduced in \cite{Lima}.  We now take the opportunity to discuss an application of the result in \cite{KSV}, which, together with \cite{KLSV}, was part of our motivation for the development of the ideas.  

In standard percolation models one associates a binary random variable (open or closed) to each vertex or bond of a given graph, and one of the basic questions regards the existence of  an infinite open connected component of containing a given vertex. A natural generalization of this consists in seeking for the existence of an infinite self-avoiding path along which one observes a given infinite binary sequence $\phi=(\phi_1,\phi_2,\dots)$. This problem was carefully set in \cite{BK}, where the authors provided, among other contributions, results for $\ZZ^d$ with large $d$. It was also studied in  \cite{KSZ1} and \cite{KSZ2} for two important classical examples: the triangular lattice and the closed packed graph of $\mathbb{Z}^2$. 

Our goal is to discuss a situation related to the kind of environment treated in this paper, and which can be seen as an application of the result in \cite{KSV}, sharing as well some of the ideas that brought to \cite{KLSV}.
%$\mathbb{L}^2=(\ZZ_+^2, \EE)$
\begin{definicao}
\label{palavra} 
Let $\Omega:=\{0,1\}^\VV$, where $\VV$ denote the set of vertices of an infinite locally finite graph. Given a sequence $\phi=(\phi_1,\phi_2,\dots) \in \{0,1\}^\NN=:\Phi $ we say that $\phi$ is seen in the configuration $\omega \in \Omega$ from the vertex $v_0$ if there exists a self-avoiding path $\gamma=(v_0,v_1, \dots)$ in the graph (i.e. the $v_i$ are all distinct and for each $i$, $v_i$ and $v_{i+1}$ are adjacent) such that $\omega(v_i)=\phi_i$ for all $i \ge 1$. One says that $\phi$ is seen in $\omega$ if there exists a vertex $v_0$ such that $\xi$ is seen in $\omega$ from $v_0$.
\end{definicao}

%$\Omega$ is thought as the environmental space while $\Phi$ is the space of words.

For $\alpha, \beta \in [0,1]$, we now consider the corresponding Bernoulli measures $\PP_{\alpha}$ on $(\Omega,\A)$ (the environment)  and $\mu_\beta$ on $(\Phi,\B)$ (space of words), where $\A$ and $\B$ denote the corresponding Borel $\sigma$-fields. That is, under $\PP_\alpha$, the $\omega (v), v \in \VV$ are i.i.d. with $\PP_{\alpha}(\omega(v)=1)=\alpha$. Similarly, under $\mu_\beta$, the $\phi_i, i \in \NN$ are i.i.d. with $\mu_\beta(\phi_i=1)=\beta$.

\medskip

As proven in \cite{BK}, one knows that 

\noindent $\bullet$ $\Lambda:=\{(\phi,\omega)\colon \phi \text{ is seen in }\omega\} \in \B\times \A$. %%\marginpar{\tiny ME Preciso verificar o quão geral são os grafos em \cite{BK}.\blue{Basta localmente finito}}

\noindent $\bullet$ Let $g_\alpha (\phi):=\PP_\alpha\{\omega \colon \phi \text{ is seen in }\omega\}$.
Then $g_\alpha(\phi) \in \{0,1\}$, for all $\phi \in\Phi$. When $g_\alpha(\phi)=1$ we just say that {\it $\phi$ percolates}.

\noindent $\bullet$ Let $h(\alpha,\beta):=\mu_\beta\{\phi \in \Phi\colon g_\alpha(\phi)=1\}$. 
Then $h(\alpha,\beta)\in \{0,1\}$. When $h(\alpha,\beta)=1$ we simply say that {\it the random word percolates}.

Given a concrete graph, one might hope to study the function $h$. Here we focus on the graph $\ZZ_+^2$
mentioned in the beginning of the paper, to which we now give an orientation for the edges, i.e. the only allowed paths are those that move at each step one unit upwards or one unit to the right. If we were thinking of the word $\phi=(1,1,\dots)$ the its percolation would correspond to the usual oriented site percolation on $\ZZ_+^2$, for which one knows there is a critical value ${\tilde p}_c \in (0,1)$ %%\marginpar{\tiny ME acertar notação}

From the above definitions, the following general properties follow at once: 

\begin{proposition}
\label{palavra-remark} Let $G$ a locally finite graph and $p_c(G)=\sup\{\alpha\in[0,1] \colon\\ g_\alpha (1,1,\dots)=0\}$, the site percolation threshold of $G$. Hence, it holds: 

\noindent (i) $h(\alpha,\beta)=h(1-\alpha,1-\beta)$ for all $(\alpha,\beta) \in [0,1]^2$.

\noindent (ii) $$h(\alpha, 1)= 
\begin{cases} 1, &\text{if } \; \alpha>p_c(G), \\
0, \; &\text{if } \; \alpha < p_c(G).
\end{cases}
$$

\noindent (iii) 
$$h(1, \beta)= 
\begin{cases} 1, &\text{if } \; \beta=1, \\
0, \; &\text{if } \;\beta \neq 1.
\end{cases}
$$

\noindent (iv) If $G=\ZZ_+^d$, with $d \ge 2$, then
\begin{equation}
\label {palavra-eq1}
    h(\alpha,\beta)=0 \text{ if  } \alpha \beta+ (1-\alpha)(1-\beta) \leq \frac{1}{d}
\end{equation}
\end{proposition}

\begin{proof} 

Itens (i)-(iii) are trivial. For the proof of (iv), let us consider the product space $\Omega \times \Phi, \A \times \B, \PP_\alpha \times \mu_\beta)$ where we take the filtration ${({\F}_n)}_n$, with ${\F}_n$ the $\sigma$-field determined by the variables $(\phi_i, i \leq n)$ and the $(\omega(x), \|x\|_1\leq n)$, where $\|\cdot\|_1$ is the  $\ell_1$-distance on $\ZZ^d$. If $Z_n$ denotes the cardinality of the set of vertices $v$ with $\|v\|_1=n$ and for which there is a path $\gamma=(0,v_1, \dots, v_n)$ with $v_n=v$ and such that $\omega(v_i)=\phi_i$ for all $i=1, \dots, n$. Then $Z_n$ is clearly ${\F}_n$-measurable and integrable. On the other hand, a simple computation gives that (writing $E$ for the expectation with respect to $\PP_\alpha \times \mu_\beta$)
$$
E(Z_{n+1}|{\F}_n, \phi_{n+1})=  (1_{[\phi_{n+1}=0]} d(1-\alpha) +1_{[\phi_{n+1}=1]} d\alpha) Z_n
$$
from which it follows that $E(Z_{n+1}| {\F}_{n})\le d c Z_{n}$, so that under the condition on (iv), the sequence $\{Z_n\}_n$ is a positive supermartingale with $E(Z_1) \leq 1$. This implies that $Z_n$ tends a.s. to 0 under $\PP_\alpha \times \mu_\beta$. Thus  
$$
\PP_\alpha \times \mu_\beta\{(\omega, \phi) \colon \phi \text{ is seen in }\omega\}=0.
$$
We now use Fubini's theorem to write
\begin{eqnarray*}
0&=&\PP_\alpha \times\mu_\beta \{(\omega,\phi)\colon \phi \text{ is seen in }\omega\}\\
&=&\int_{\Phi}\int_{\Omega} 1_{\{(\omega, \phi) \colon \phi \text{ is seen in } \omega\}}d\PP_\alpha d\mu_\beta \\
&=& \int_\Phi \PP_\alpha \{\omega \colon \phi \text { is seen in } \omega \} d\mu_\beta\\
&=&\int_\Phi g_\alpha (\phi) d\mu_\beta(\phi)=\mu_\beta\{\phi \colon  g_\alpha(\phi)=1\}=h(\alpha,\beta),
\end{eqnarray*}
which concludes the proof.
\end{proof} 

Going to the specific case of $Z^2_+$ with all the edges being oriented in the increasing direction,  one can add the following result as a direct consequence of Theorem 1.1 in \cite{KSV}. 

\begin{proposition}
\label{palavra-prop}
 Let $G=\ZZ_+^2$, oriented as described above. If $\alpha > \tilde{p}_c$ one can find $\beta_0(\alpha)<1$ such that $h(\alpha,\beta)=1$ for all $\beta\geq \beta_0(\alpha)$.
\end{proposition}

\begin{proof}
Given the pair $(\omega, \phi)$, we say that $v$ is open is $\omega (v)=\phi_{\|v\|_1}$ and closed otherwise. We also declare each line $r_n=\{v \in \ZZ_+^2\colon \|v\|=n\}$ as \emph{ good }or \emph{bad} according to $\phi_n=1$ or $\phi_n=0$.\footnote{The attribute good/bad here is exactly the opposite of that used for the layers in Section \ref{intro}.} It is easy to see that this defines an oriented site percolation model as investigated in \cite{KSV}. We have 
$$
\delta=P\{r \text{ is bad }\}=1-\beta
$$
$$
p_g=P\{ v \text{ is open }| r_{\|v\|_1} \text{ is good }\}=\alpha 
$$
$$
p_b=P\{ v \text{ is open }| r_{\|v\|_1} \text{ is bad }\}=1-\alpha
$$
On the other hand, Theorem 1.1 from \cite{KSV} tells us that for each $\alpha >\tilde p_c$, there exists 
$\beta_0(\alpha)< 1$ such that for all $\beta>\beta_0$ there is a positive probability of an infinite open path starting at the origin. According to our definitions, this implies 
$$
\PP_\alpha \times \mu_\beta \{(\omega, \phi) \colon \phi \text{ is seen in } w\}=1,  \text { for all } \beta > \beta_0(\alpha),
$$
for all $\beta>\beta_0(\alpha)$. From this, and using Fubini's theorem as above we conclude the proof.
\end{proof} 

\begin{remark}
\label{final-ksz}

\noindent (i) Of course, the statement of Proposition \ref{palavra-prop} extends to $\ZZ_+^d$ (oriented) for all $d \ge 2$. 

\noindent (ii) In \cite{KSZ1}, the authors considered the triangular lattice  (which is $\ZZ^2$ with added edges between $(i,j)$ and $(i+1,j+1)$ for each $i,j$. In this case the critical parameter for site percolation is $1/2$. They proved that for $\alpha=1/2$ and $\beta \in (0,1)$, $\mu_\beta$ almost all words are seen. 

\noindent (iii) In \cite{KSZ2} the authors studied percolation of random words on the graph  $\ZZ^2_{cp}$ (usually called closed packed $\ZZ^2)$, i.e. there are edges added between any pair $(i,j)$ and $(i-1,j+1)$, as well between  $(i,j)$ and $(i+1,j+1)$. This is a graph for which the critical parameter for site percolation is less than 1/2. In particular, when $\alpha=1/2$ the words $(1,1,\dots)$ and $(0,0,\dots)$ percolate. In \cite{KSZ1} they proved that if $p_c(\ZZ^2_{cp}) < \alpha < 1-p_c(\ZZ^2_{cp})$ than all the words are seen, i.e.
$$
\PP_\alpha\{ \omega \colon \text{all words are seen in } \omega\}=1
$$
It makes sense to investigate whether something similar holds for other graphs for which the critical parameter for site percolation is smaller than 1/2. As mentioned earlier, results for large dimensions were given in \cite{BK} and recently extended for $\ZZ^d$ for all $d\geq 3$ in \cite{NTT}.
\end{remark}

\noindent {\bf Acknowledgements.}

The research that brought to this paper was initiated during the preparation of the Ph.D. thesis of B.N.B.L. It continued over several years that included many inspiring discussions with Harry Kesten, to whom we are deeply indebted. 

% from discussions during various meetings at CBBF, IMPA, UFMG, and UFRJ, as well as visits of M.E.V. and V.S. to MSRI. B.N.B.L. and M.E.V. thank all these institutions for the hospitality and support. 

During the initial period, B.N.B.L. was partially supported by CNPq grant 301844/2008-9, and M.E.V. was partially supported by CNPq grant 302796/2002-9. During the current revision process,  B.N.B.L. is partially supported by CNPq grant 305811/2018-5 and FAPEMIG (Programa Pesquisador Mineiro)  and M.E.V. is partially supported by CNPq grant 310734/2021-5 and Faperj grant E-26/202.636/2019.


\bigskip


\begin{thebibliography}{99}

\bibitem{ACCN} Aizenman, M., Chayes, J.T., Chayes, L., and Newman, C.M. (1988) Discontinuity of the magnetization in one-dimendional $\frac{1}{{|x-y|}^2}$ Ising and Potts models. {\em Journal of Statistical Physics} {\bf 50}, 1-40.

%\bibitem{BBS}Balister, P.N., Bollobas, B., and Stacey, A.M. (2000) Dependent Percolation in two dimensions. {\em Probability Theory and Related Fields} {\bf 117}, 495-513.

\bibitem{BK}Benjamini, I. and Kesten, H. (1995) Percolation of arbitrary words in $\{0,1\}^\NN$. {\em Annals of Probability} {\bf 23}, 1024-1060.
%\bibitem{BPP}Benjamini, I., Pemantle, R., and Peres, Y. (1998) Unpredictable paths and percolation. {\em Annals of Probability} {\bf 26}, 1198-1211.


\bibitem{BDS}Bramson, M., Durrett, R.T., and Schonmann, R.H. (1991) The contact process in a random environment. {\em Annals of Probability} {\bf 19}, 960-983.

\bibitem{BS} Basu, R. and Sly, A. (2014) Lipschitz embeddings of random sequences. {\em Probability Theory and Related Fields} {\bf 159}, 721-775.

\bibitem{CS} Chayes, L. and Schonmann, R.H. (2000) Mixed percolation as a bridge between site and bond percolation. {\em The Annals of Applied Probability} {\bf 10}, 1182-1196.

%\bibitem{CD} Cox, J.T. and Durrett, R.T. (1983) Oriented percolation in dimensions $d\geq 4$: bounds and asymptotic formulas. {\em Mathmatical Proceedings of the Cambridge Philosophical Society} {\bf 93}, 151-162.

\bibitem{DHKS} Duminil-Copin, H., Hil\'ario, M.R., Kozma, G., and Sidoravicius, V. (2018) Brochette percolation. {\em Israel Journal of Mathematics} {\bf 225}, 479-501.

%\bibitem{Du} Durret, R. (1984) Oriented percolation in two dimensions. {\em The Annals of Probability} {\bf 12}, 999-1040.

\bibitem{Gr}Grimmett, G.R. (1999) Percolation, 2nd edition. Springer-Verlag, Berlin.

\bibitem{GHM} Georgii, H.O., H\"aggstr\"om, O., and Maes, C. (2001) The random geometry of the equilibrium phases In {\em Phase Transition and Critical Phenomena} (C. Domb and J. Lebowitz, ed.) {\bf 18} 1-142. Academic Press, London.

\bibitem{Ho}Hoffman, C. (2005) Phase transition in dependent percolation. {\em  Communications in Mathmatical Physics} {\bf 254}, no 1  1-22.

\bibitem{JMP}Jonasson, J., Mossel, E., and Peres, Y. (2000) Percolation in dependent random environment. {\em Random Structures and Algorithms} {\bf 16}, 333-343.
%\bibitem{Ke1}Kesten, H. (1980) The critical probability of bond percolation on the square lattice equals $\frac{1}{2}$.
%{\em Communications in Mathmatical Physics} {\bf 74}, 41-59.

%\bibitem{Ke2}Kesten, H. (1982) Percolation theory for Mathmaticians. Birkh\"auser, Boston.

\bibitem{KLSV}Kesten, H., de Lima, B.N.B., Sidoravicius, V., and Vares, M.E. (2014) On the compatibility of  binary sequences. {\em Communications on Pure and Applied Mathematics} {\bf 67}, 871-905.

\bibitem{KSV}Kesten, H., Sidoravicius, V., and Vares, M.E. (2022) Oriented percolation in a random environment. To appear in {\em Electronic Journal of Probability}.

\bibitem{KSZ1} Kesten, H., Sidoravicius, V. and Zhang, Y. (1998)
    Almost all words are seen in critical site percolation on the triangular lattice.
   {\em  Electronic Journal of Probability} {\bf 3} n. 10, 75 pp.

\bibitem{KSZ2}Kesten, H., Sidoravicius, V., and Zhang, Y. (2001) Percolation of arbitrary words on the closed packed graph of $\ZZ^2$. {\em Eletronic Journal of Probability} {\bf 6} n. 4, 1-27.

\bibitem{Kl}Klein, A. (1994) Multiscale Analysis in disorder systems: percolation and contact process in random environments. In {\em Probability and Phase transitions}, Nato ASI series, Kluwer Academic Publishers, Dordecht.

%\bibitem{Li2}Liggett, T. (1985) Interacting Particle Systems. Springer-Verlag, Berlin.

\bibitem{Li}Liggett, T. (1992) The survival of the one dimensional contact process in random environments. {\em Annals of Probability} {\bf 20}, 696-723.

\bibitem{Lima} de Lima, B.N.B.  {\em Percolação de Bernoulli dependente em $\mathbb{Z}^2$.} Ph.D.Thesis, 2003. Informes de Matemática. IMPA, Série C-26/2004.
% \bibitem{MW4} McCoy, B.M. (1969) Incompleteness of the Critical Exponent Description for Ferromagnetic Systems Containing Random Impurities. {\em Physical Review Letters} {\bf 23} n. 7, 383-386.
\bibitem{MW3} McCoy, B.M. (1969) Theory of a Two-Dimensional Ising Model with Random Impurities. III. Boundary Efects. {\em Physical Review} {\bf 188} n. 2, 1014--1031.

\bibitem{MW4} McCoy, B.M. (1970) {\em Theory of a two-dimensional Ising model with random impurities. IV. Generalizations}. Phys. Rev. B {\bf 2}, 2795--2803.


\bibitem{MW1} McCoy, B.M. and Wu, T.T. (1968) Theory of a Two-Dimensional Ising Model with Random Impurities. I. Thermodynamics. {\em Physical Review} {\bf 176} n. 2, 631-643.

\bibitem{MW2} McCoy, B.M. and Wu, T.T. (1969) Theory of a Two-Dimensional Ising Model with Random Impurities. II. Spin Correlation Functions. {\em Physical Review} {\bf 188} n. 2, 982-1013.

\bibitem{NTT} Nolin, P., Tassion, V., and Teixeira, A. (2022) No exceptional words for Bernoulli percolation. To appear in {\em Journal of the European Mathematical Society}.

\bibitem{Wi} Winkler, P. (2000) Dependent percolation and colliding random walks. {\em Random Structures and Algorithms} {\bf 16}, 58-84.
\end{thebibliography}
\end{document}